%% file: supermarketoptimization.tex
\documentclass[]{article}

\usepackage{graphicx}
\usepackage{subcaption}

\usepackage{setspace}
\usepackage{tabularx}
\usepackage{longtable}
\usepackage{booktabs}
\usepackage{array}

\usepackage{graphicx}
\usepackage{enumitem}
\usepackage{multirow}
\usepackage{amsmath}
\usepackage{amssymb}
\usepackage{moreverb,url}
\usepackage[margin=1in,headsep=.3in]{geometry}
\usepackage{csquotes}
\usepackage{longtable}
\usepackage{authblk}
\usepackage{cite}
\usepackage{hyperref} \hypersetup{colorlinks}

\title{Hundreds of grocery outlets needed across the United States to achieve walkable cities}

\author[1]{D Horton}
\author[2]{T M Logan}
\author[1]{E Speakman}
\author[3]{D Skipper}

\affil[1]{Mathematical and Statistical Sciences, University of Colorado, Denver}
\affil[2]{Civil and Natural Resources Engineering, University of Canterbury, New Zealand}
\affil[3]{Mathematics Department, United States Naval Academy, Annapolis, Maryland}

\date{}

\begin{document}
\maketitle

\section*{Abstract}
The notion of the $x$-minute city is again popular in urban planning, but the practical implications of developing walkable neighborhoods have not been rigorously explored.
What is the scale of the challenge that cities needing to retrofit face?
Where should new stores or amenities be located?
For 500 cities in the United States, we explored how many additional supermarkets would be required to achieve various levels of $x$-minute access and where new stores should be located so that this access is equally-distributed.
Our method is unique because it combines a novel measure of equality with a new model that optimally locates amenities for inequality-minimizing community access.
We found that 25\% of the studied cities could reach 15-minute access by adding five or fewer stores, while only 10\% of the cities could even achieve 5-minute average access when using neighborhood centroids as potential sites; the cities that could, on average, required more than 100 stores each.
This work provides a tool for cities to use evidenced-based planning to efficiently retrofit in order to enable active transport, benefiting both the climate and their residents' health.
It also highlights the major challenge facing our cities due to the existing and ongoing car-dependent urban design that renders these goals unfeasible.

\bigskip

\section{Introduction}
The presence of amenities in urban areas is a key enabler of active transport and sustainable urban design \cite{McCormack2008-nl}.
This type of urban design has benefits including public health \cite{Celis-Morales2017-mn, Lindsay2011-ys, Stevenson2016-yb, Mueller2020-bl, Staricco2022-zk}, improved quality of life \cite{Leyden2003-dg, Lopez2011-cc}, emissions reductions \cite{Brand2021-qw, Iea2020-is}, and resilience through improved social cohesion \cite{Dempsey2011-og, Jacobs1961-po, Hao2022-tp, Logan2021-inequ, Anderson2022-vr, Sajjad2021-vq}.
However, many cities around the world do not have a sufficient density of amenities to promote active transport (walking and cycling), as their urban design has largely been car-oriented \cite{Logan2022-mr, Wu2021-nx,Pucher2010-ip}.
Between 1990 and 2014, urban sprawl increased globally by 95\% \cite{Behnisch2022-xa}.
Cities worldwide are now articulating visions (e.g., the 10-minute city) of improving their residents' proximity to amenities and thereby capitalizing on all of the benefits of active transport \cite{C40_Cities2020-bm}.
This sustainable transition raises significant questions for our existing cities about how they can direct this retrofit efficiently and effectively.

Simultaneously, cities are generally working to achieve this transition in a manner that addresses and mitigates inequities.
The urban design of a community/city influences how resources and burdens are distributed between residents \cite{Marino2012-lq,Wilson2008-yk, Calvin2017-ja}.
Unfortunately, empirical evidence indicates that these burdens and resources are currently not equally distributed among people;
globally, disadvantaged and underprivileged people are systematically exposed to larger environmental burdens and have lower access to beneficial resources \cite{Fussel2010-te, Bulkeley2014-so, White2021-jd, Logan2021-inequ}.
This is a form of distributive injustice \cite{Low2013-yx}.
A novel approach to measuring equality of access and equity of access between different groups (e.g., socioeconomic, demographic, etc.) was introduced in 2020: the Kolm-Pollak equally-distributed equivalent (EDE), a measure similar to the Atkinson Index (commonly used to evaluate income inequality between countries) that is suitable for urban contexts \cite{Sheriff2020-ge, Logan2021-inequ}.
The EDE is essentially an inequality-penalized average/mean for a statistical distribution.
That is, it can be used to express the average distance of a community, but with penalties on distances that are higher than the mean  (i.e., the tail of the distribution) so that it better represents the actual experience of the residents.

Using this metric, we evaluated access to grocery stores and supermarkets in the 500 largest cities in the United States (US). \autoref{tab:cities_assess} shows the largest 20 cities in the study, along with their residents' inequality-penalized average access, and their ``access to supermarkets'' ranking among the 500 cities we assessed.
Note that we use the words grocery stores and supermarkets interchangeably; by both we mean stores that sell food, including fresh produce, and are larger than convenience stores or gas stations.

\begin{table}[ht]
\setlength{\tabcolsep}{7pt}
\centering
    \caption{
       Summary results for the 20 largest US cities. ``Rank'' indicates the rank of the city with respect to supermarket access among the 500 largest cities in the US (1 is best). ``EDE Distance (km)'' indicates the equity-penalized mean distance (EDE) of residents to a grocery store. See
       the supplementary materials for the full ranked list of 500 US cities.
    }
    \begin{tabular}{lccr}
    \toprule
        \textbf{City, State} & \textbf{Rank} & \textbf{EDE Distance (km)} & \multicolumn{1}{c}{\textbf{Metro Population}} \\
        \midrule
        New York, NY      & 3   & 0.8 & 8,784,592 \\
        San Francisco, CA & 7   & 1.0 & 871,136   \\
        Philadelphia, PA  & 15  & 1.1 & 1,593,147 \\
        Washington, DC    & 17  & 1.1 & 684,900   \\
        Chicago, IL       & 20  & 1.1 & 2,733,239 \\
        \midrule
        Seattle, WA       & 26  & 1.1 & 726,482   \\
        Los Angeles, CA   & 56  & 1.3 & 3,849,235 \\
        San Jose, CA      & 99  & 1.6 & 993,779   \\
        Denver, CO        & 110 & 1.6 & 705,515   \\
        San Diego, CA     & 143 & 1.7 & 1,347,374 \\
        \midrule
        Houston, TX       & 209 & 1.9 & 2,215,641 \\
        Charlotte, NC     & 232 & 2.0 & 804,437   \\
        Columbus, OH      & 260 & 2.1 & 868,417   \\
        Dallas, TX        & 275 & 2.1 & 1,269,024 \\
        Phoenix, AZ       & 321 & 2.3 & 1,553,053 \\
        \midrule
        Indianapolis, IN  & 370 & 2.5 & 788,869   \\
        San Antonio, TX   & 388 & 2.6 & 1,381,080 \\
        Fort Worth, TX    & 449 & 3.3 & 865,707   \\
        Jacksonville, FL  & 483 & 4.1 & 834,225   \\
        Austin, TX        & 489 & 4.8 & 893,947   \\
        \bottomrule
    \end{tabular}

    \label{tab:cities_assess}
\end{table}

While this ranking and evaluation is interesting, the key question for planners and urban and justice advocates is how can access in existing urban areas be improved in an effective and equitable manner?
We address this question by asking two specific questions for each of the 500 largest cities in the US:
\begin{enumerate}[label=(\arabic*)]
        \item If a city can open $k$ additional supermarkets, where should they be located to best improve equitable access?
        \item If we want to reach some level of equitable access (e.g., 15 minutes), how many additional supermarkets are required, and where should they be located?
    \end{enumerate}
The answers to these questions provide an indication of the scale of change that is required to retrofit car-oriented cities to enable active transport modes and move towards health-promoting and sustainable urban design.

We address these questions by developing an optimization approach based on the measure of inequality, the EDE.
The technical mathematical advances are detailed in a sister article \cite{horton2024scalable}.
In this paper, we apply those methods to explore the scale of the retrofit required for the 500 largest cities in the US, considering grocery stores.

While this optimization approach is general to any amenity or destination type, we selected grocery stores because they are a commonly frequented destination and because of the prevalence of food deserts worldwide.
A food desert is defined in the US as a region more than 1 mile (1600 meters) away from the nearest grocery store (in an urban area) \cite{foodinescurity}.
Food deserts contribute to food insecurity, which is the state of being without reliable access to a sufficient quantity of affordable, nutritious food \cite{foodinescurity}.
Healthy food access is a factor in mitigating chronic disease \cite{Kolak2018-az, Garcia2020-xt,Apparicio2007-di, Walker2010-ch} and is an environmental justice issue due to the disproportionate impacts on racial/ethnic minority and low-income communities \cite{Krenichyn2006-ve,Day2006-ak, Kolak2018-az, Walker2010-ch}.
The Covid-19 pandemic compounded this issue by increasing the prevalence of food insecurity \cite{MORALES2021} and research indicates similar socio-demographic determinants of food insecurity and infection rates, notably among Black and American Indian populations \cite{Kimani2021-xx}.
Additionally, in low-income households with children, there was a 22\% increase in food insecurity from 2019 to 2020 \cite{Sharma2020-ol}.

Traditional optimization models, including those designed to optimally locate amenities for residential access, were originally developed with commercial applications in mind \cite{pmed1, pmed2, covering}.  As such, they focus on minimizing costs or maximizing profits and do not consider equitable access for the population. For example, facility location optimization models often minimize the mean distance or travel time between facilities and demand points, resulting in an overall reduction of transportation costs \cite{pmed1}. The drawback of the mean-minimizing model in an equity context is that it sometimes leads to solutions where the minimum average distance is attained by improving access a little for many people (by placing more stores in heavily populated areas that already have supermarkets), rather than targeting those who are currently disadvantaged. In recent years, there has been significant work aimed at incorporating equity into facility optimization models \cite{Karsu2015-cb}.  Unfortunately, equity metrics tend to be algebraically complex, so optimization models that contain them do not scale computationally to practical problem sizes \cite{Barbati2018-nd}.
The model detailed in the sister article \cite{horton2024scalable} overcomes this limitation, enabling us to apply the model to the 500 cities in our study.

\section{Methods}
In this paper, we seek to evaluate and optimize access (and access inequality) to grocery stores across the 500 largest cities in the US.
In this section, we describe our data, provide background information on the metric we use to quantify inequality-penalized access, and present our optimization models.

\subsection{Cities}

We selected the 500 largest cities in the US based on 2020 US Census population data \cite{census2022}.

\subsection{Measuring access}

We calculate the driving and walking distance from the centroid of the US Census Block (the smallest census unit for the US) to all existing grocery stores and potential store locations.

We measure access as the distance to the nearest amenity.
This means that our approach is not considering the demand for a particular amenity or it's capacity to serve that demand and is a limitation that we will seek to address in future.
However, while not the only requirement, proximity to services is necessary for access \cite{Penchansky1981-qh, Saurman2016-gj}.

To calculate the distance to the nearest amenity, we utilize the method described by Logan et al. \cite{Logan2019-fr}.
This leverages the Open Source Routing Machine \cite{osrm} to calculate the network distances between origins (Census Blocks) and destinations (existing and potential store locations).
This method accounts for geographical barriers, such as freeways, waterways, and railroad tracks.

This method does not account for the suitability or quality of the walking environment (see \cite{Koschinsky2016-ao} for a discussion) but only whether a route exists.

\textit{Grocery store locations.}
We use existing supermarket locations within a 5km radius of the city from the USDA's Food and Nutrient Service SNAP database available on ArcGIS Hub.
This is consistent with the analysis of \cite{Logan2022-mr}.

\textit{Potential grocery store locations}.
We used centroids of US Census Block Groups from the 2020 US Census to geographically cover each city with potential store locations. After Blocks, Block Groups are the second most granular geographic unit captured in the US Census.

\textit{Population data.}
The population of each Census Block was based on the 2020 US Census and exported from the IPUMS National Historical Geographic Information System \cite{nhgis}.

\subsection{Inequality metric}

The environmental justice (EJ) community has focused recently on ranking distributions of disamenities, such as pollution exposure, with the goal of quantifying and comparing the health risks that communities face. Equally distributed equivalents seek to answer the question, ``what level of risk would make an individual indifferent between a distribution in which everyone receives that risk and the actual unequal risk distribution?''
The Kolm-Pollak equally distributed equivalent (EDE) was introduced as the only metric that satisfies several key properties of ranking functions identified by the EJ community \cite{Sheriff2020-ge}. The EDE incorporates inequality by measuring the center of the distribution with a penalty for values that are above (worse) than the mean. In this way, the EDE is a more accurate measures of the actual experience of a population than the population mean. (Consider, for example, how the mean of a distribution of incomes can be very misleading.) A recent article presents a case study of 10 US cities to demonstrate how the EDE can be used to rank cities with respect to access to amenities and shows how the rankings change with the level of aversion to inequality \cite{Logan2021-inequ}.

Like other equally distributed equivalents, the Kolm-Pollak EDE depends on a user-defined parameter, $\epsilon \in \mathbb{R}$.
 If larger values in the distribution are undesirable, such as pollution level or distance to a grocery store, then $\epsilon < 0$ and the EDE is always at or above the mean of the distribution. Larger values of $|\epsilon|$ represent more aversion to inequality. In typical applications, $|\epsilon|$ is assigned a value between 0.5 and 2.

We used the EDE distance (in meters) with $\epsilon = -1$ to quantify the level of access of a community to grocery stores. For a given city, let $R$ represent the set of Census Blocks and let $p_r$ represent the population of Block $r \in R$. Let $z_r$ represent the walking distance (in meters) of Block $r \in R$ to the closest grocery store. The Kolm-Pollak EDE distance of the residents of the city to supermarkets is,
\begin{align} \label{popweightedkp}
\mathcal{K}(\mathbf{z}) ~=~ -\frac{1}{\kappa}\ln\left[\frac{1}{T}\sum_{r \in R}p_re^{-\kappa z_r}\right],
\end{align}
where $\mathbf{z}$ is the vector of distances, $T:=\sum_{r \in R}p_r$ is the total population, and $\kappa := \alpha\epsilon$, where
\begin{equation} \label{eq:alpha}
    \alpha ~=~ \frac{\sum_{r \in R} p_rz_r}{\sum_{r \in R} p_rz_r^2}.
\end{equation}
The aversion to inequality, $\epsilon$, is scaled to the problem data via $\alpha$, so $\kappa$ is the appropriately-scaled aversion to inequality. We used the current distributional access to approximate the value of $\alpha$ that corresponds to the optimal distribution of distances. This allowed us to treat $\kappa$ as a parameter (constant) in our models. For a more detailed discussion of our strategy for approximating $\alpha$, please see our companion methods article \cite{horton2024scalable}.

\subsection{Optimization}

Many models aimed at incorporating equity in facility location optimization have been proposed \cite{MARSH19941,Smith2013-md, Batta2014-uj, MARSH19941, Xu2023-iy}.
Typically, these models do not scale to large problem sizes, or even to placing more than one facility. However, linear EDE-minimizing model scales to large, city-sized instances \cite{horton2024scalable}.

Adding to the notation defined above, let $S$ represent the set of existing and potential supermarket locations, and let $C \subseteq S$ represent the set of existing (current) supermarket locations. Let $d_{r,s}$ represent the walking distance (in meters) between Block $r \in R$ and location $s \in S$. Our decision variables are all binary:
\begin{align*}
    &x_s~:=~1 \text{ if a supermarket is placed at location } s \in S \text{, 0 otherwise};\\
    &y_{r,s}~:=~1 \text{ if Block } r \in R \text{ is assigned to service location } s \in S\text{, 0 otherwise}.
\end{align*}
As a function of the vector, $\mathbf{y}$, of $y_{r,s}$ variables, the EDE is,
\begin{align}
\mathcal{K}(\mathbf{y}) ~=~ -\frac{1}{\kappa}\ln\left[\frac{1}{T}\sum_{r\in R}p_re^{-\kappa \sum_{s\in S}y_{r,s}d_{r,s}}\right].
\end{align}
We can equivalently minimize the so-called linear proxy \cite{horton2024scalable},
\begin{align}
\overline{\mathcal{K}}(\mathbf{y}) ~=~ \sum_{r\in R}\sum_{s \in S} p_r y_{r,s} e^{-\kappa d_{r,s}},
\end{align}
and then convert the optimal objective value to an EDE score:
\[ \mathcal{K}(\mathbf{y}) = -\frac{1}{\kappa}\ln\left(\frac{1}{T}\overline{\mathcal{K}}(\mathbf{y})\right). \]

\subsubsection{Question 1}

To answer the question of how to optimally locate $k$ additional supermarkets, we minimize the EDE linear proxy in the objective function, \eqref{eq:kpl_objective}:

\begin{align}
    \text{minimize} \hspace{0.5cm}& \overline{\mathcal{K}}(\mathbf{y}) ~=~ \sum_{r \in R}\sum_{s \in S}p_r y_{r,s}e^{-\kappa d_{r,s}},\label{eq:kpl_objective}\\
    \text{subject to} \hspace{0.5cm}& \sum_{s \in S  \backslash C} x_s = k; \label{eq:numfacilities}\\
    & \sum_{s \in S} y_{r,s} = 1, \quad\forall ~r \in R; \label{eq:eachorgin1destination}\\
    &y_{r, s} \leq x_s, \quad\forall ~r \in R, ~s \in S; \label{eq:openfacilities}\\
    & x_s = 1, \quad\forall ~s \in C; \label{eq:existing}\\
    &\, x_s, y_{r,s} \in \{0,1\}, \quad\forall ~r \in R, ~s \in S. \label{eq:binary}
\end{align}
As noted above, we can convert the optimal objective value to an EDE distance to determine the optimal level of access that can be achieved by adding $k$ stores.
Constraint \eqref{eq:numfacilities} ensures the correct number of new supermarkets are opened.  Constraint \eqref{eq:eachorgin1destination} ensures that every Census Block is assigned to a single store, while \eqref{eq:openfacilities} ensures that the assigned store location is open. Constraint \eqref{eq:existing} keeps all existing stores open, and \eqref{eq:binary} enforces the binary requirement on the indicator variables.

\subsubsection{Question 2}

In the model that answers how many (and where) stores should be opened to achieve a given level of equitable access, the EDE linear proxy is included as a constraint.
Suppose we are aiming for a level of equitable access of no more than $\ell$ meters.
 We must convert $\ell$ to the same units as the linear proxy to serve as the upper bound on the access constraint: $L = Te^{-\kappa \ell}$. The model that answers our second question is:
\begin{align}
    \text{minimize} \hspace{0.5cm}&\sum_{s \in S\backslash C} x_s, \label{eq:min_stores}\\
    \text{subject to}\hspace{0.5cm}&\sum_{r \in R} \sum_{s \in S} p_r y_{r, s}e^{-\kappa d_{r, s}} \leq L; \label{eq:meet_access}\\
    & \sum_{s \in S} y_{r,s} = 1,  \quad\forall ~r \in R; \tag{\ref{eq:eachorgin1destination}}\\
    &y_{r, s} \leq x_s, \quad\forall ~r \in R, ~s \in S; \tag{\ref{eq:openfacilities}}\\
    & x_s = 1, \quad\forall~s \in C; \tag{\ref{eq:existing}}\\
    &\, x_s, y_{r,s} \in \{0,1\},  \quad\forall \, r \in R, ~s \in S. \tag{\ref{eq:binary}}
\end{align}
The objective function, \eqref{eq:min_stores}, minimizes the number of new stores, while \eqref{eq:meet_access} ensures the desired level of access is achieved. The rest of the constraints are the same as in the previous model.

\subsection{Computing environment}

We implemented the models in Python using the optimization modeling language Pyomo \cite{hart2011pyomo,bynum2021pyomo}, and solved the models using the linear mixed-integer optimization solver, Gurobi \cite{gurobi}. We solved most instances on a high-performance computing cluster, an Advanced Micro Devices (AMD) 7502 CPU processor with 64 cores and 512 GB of memory, allocating one out of the 64 available cores to each instance. The New York instances required more memory. For those, we used an AMD 7502 CPU processor with 64 cores and 2 TB of memory.

\section{Results}

\subsection{Where should new facilities be located?}
The first question we address is: if a city can open $k$ additional amenities, where should they be located to best improve equitable access?
To answer this question, we developed an approach that minimizes the EDE: a metric that captures the average of a quantity (such as minimizing the distance to nearest supermarket) but penalizes for inequality \cite{horton2024scalable}.
In this study, we applied this method to the 500 largest cities in the US.

\begin{figure}
    \centering
    \includegraphics[width=\textwidth]{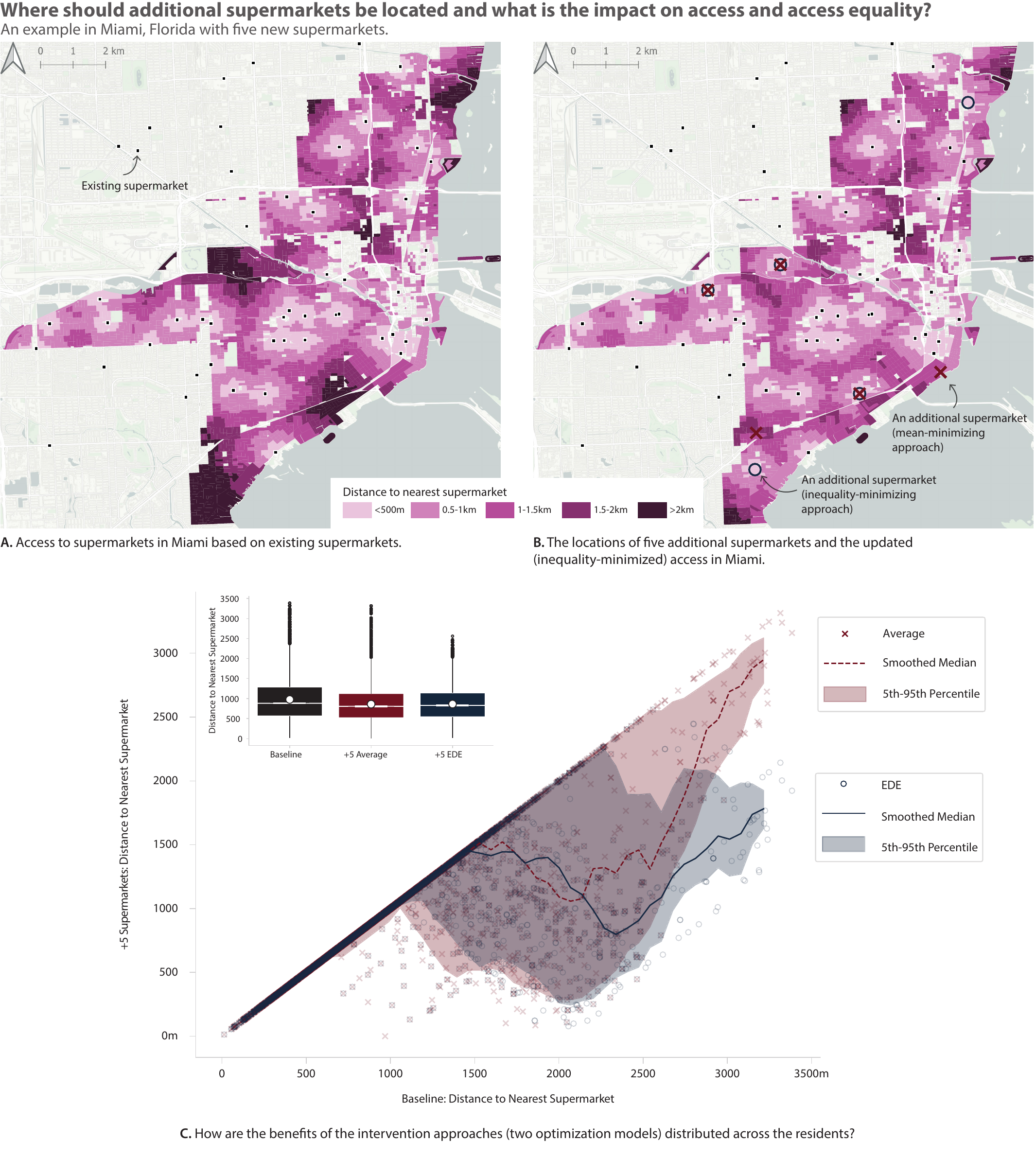}
    \caption{Where should additional supermarkets be located to improve access and access equality? Maps A and B show the access to supermarkets in Miami before and after the addition of five equality-optimizing supermarkets. Map B also shows the mean-minimizing locations of five additional supermarkets. The graphics in C show that the EDE-minimizing approach best targets the residents who currently have the worst access. The main graphic in C shows which Blocks benefited from each intervention, while the boxplot shows population weighted means (indicated by a circle) and quartiles of before and after each intervention.}
    \label{fig:Miami}
\end{figure}

For example, \autoref{fig:Miami}A shows a map of Miami, Florida, with shading to indicate the distance of residents to supermarkets based on the locations of supermarkets at the date of this study.
In \autoref{fig:Miami}B, we show the recommended locations for five new supermarkets based on the traditional (mean distance minimizing) and our proposed equitable optimization approaches. \autoref{fig:Miami}B is shaded according to the updated distance to nearest supermarket given the optimal locations proposed by our method. Even though some of the added stores are sited at or near the same location under both approaches, our proposed placement leads to notable improvements for equality vs the traditional approach.

In order to visualize the distributional effect of each intervention, we plot the access for each Block before and after the optimally-located supermarkets in \autoref{fig:Miami}C. The main graphic includes two marks for each Census Block in Miami: a red ``$\times$'' corresponding to the mean-minimizing approach and a blue
  ``\begingroup\makeatletter\def\f@size{9}\check@mathfonts\raisebox{\depth}{$\boldsymbol{\bigcirc}$}\endgroup'' corresponding to the inequality-optimizing approach.
When a point is shown on the 1:1 line it means that the Block's access has not changed as a result of the intervention.
If a point is above the 1:1 line it means that a Block's access has become worse due to the intervention, whereas if it is below the 1:1 line, the access has improved.
This figure shows the distribution effect because it shows which Blocks experience the greatest improvement in access under each intervention.
By comparing the traditional vs inequality-minimizing approaches, we see that the approach that seeks to optimize inequality tends to improve the access of currently access-poor areas in comparison to the traditional approach.
This is because the mean/average can be minimized by improving the access (reducing the distance) of any area.

The box-plot in the upper left corner of \autoref{fig:Miami}C provides a visual representation of access statistics before and after each intervention. The inequality-minimizing (EDE-optimizing) method achieves nearly the same average and median access as the traditional (mean-minimizing) method, while more successfully targeting those with the poorest access in the baseline distribution. This effect is further analyzed and verified in the sister methods article \cite{horton2024scalable}.

 \begin{figure}
    \centering
    \includegraphics[width=0.5\textwidth]{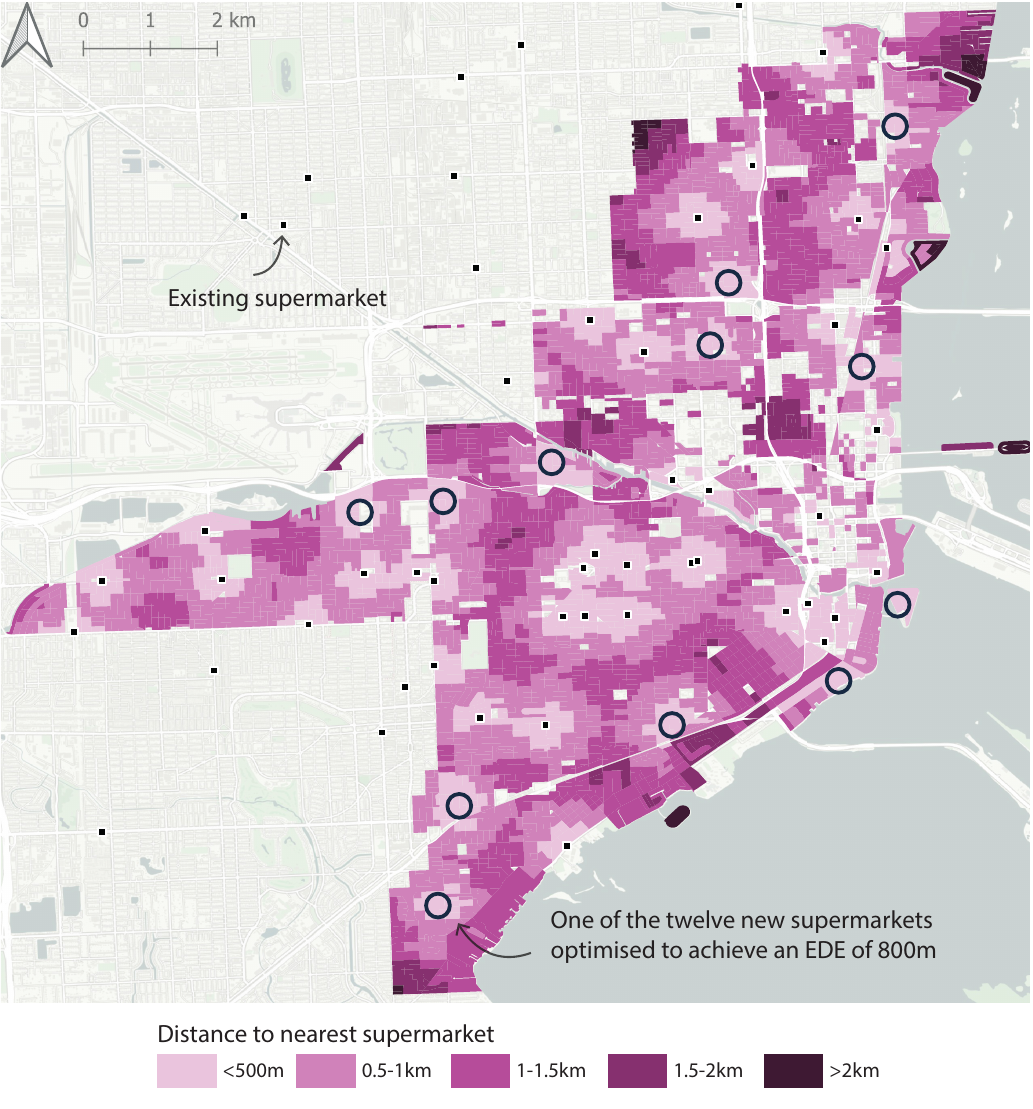}
    \caption{An equitable 10-minute plan for Miami.}
    \label{fig:miami-10min}
\end{figure}

\begin{figure}
    \centering
    \includegraphics[width=0.9\textwidth]{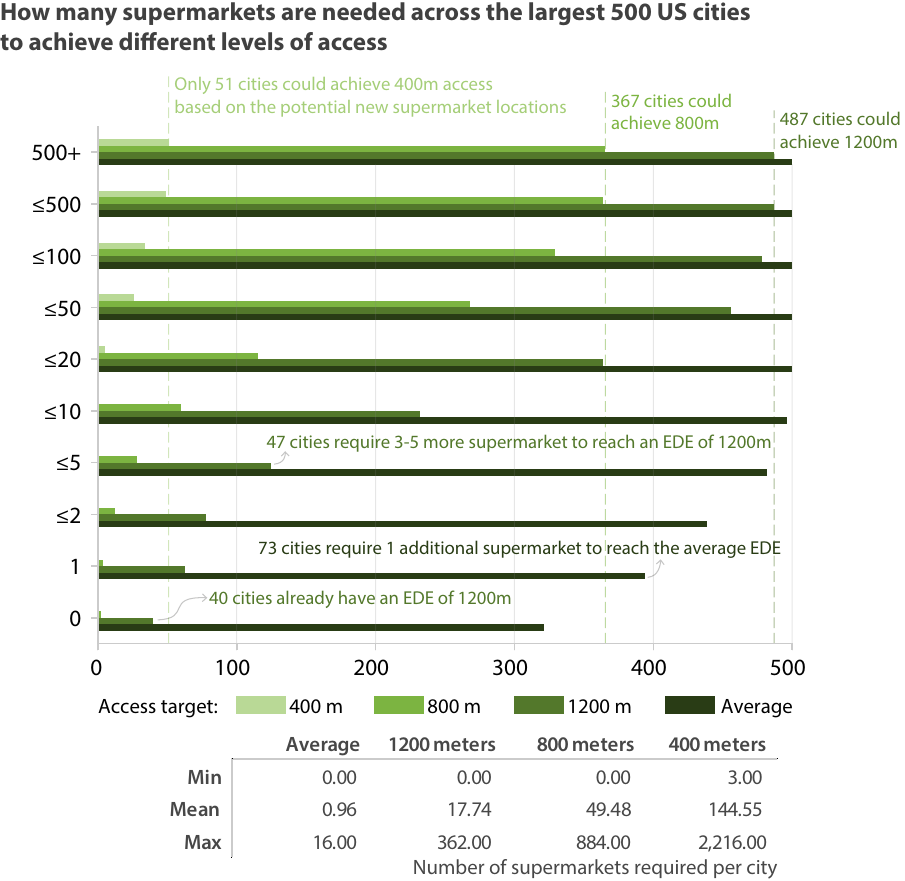}
    \caption{Nearly every city can reach the average level of access by adding fewer than 10 supermarkets. The outlook is less positive for more ambitious access levels. The statistics in the table represent cities for which the target access is feasible.}
    \label{fig:number-required}
\end{figure}

\subsection{How many facilities are required?}
Initially, we sought to answer the question of \textit{where} a city should build additional amenities to improve equitable access.
But this raises the question of \textit{how many} are required to provide a certain level of equitable access?
For instance, if a city is aiming for 10-minute walkable neighbourhoods, enabling people to reach the identified amenities within a 10-minute walk of their residence, how many amenities are required and where should they be built?

In the case of Miami, Florida, the city's supermarket access EDE is 1040 meters, which is roughly a 12-13 minute walk.
However, for Miami to decrease this to 10-minutes (an EDE of 800 meters), they would need an additional 12 stores, the locations of which are shown in \autoref{fig:miami-10min}.
If they want to decrease the travel time to 5-minutes (400 meters), they would need more than 100 additional supermarkets.

However, Miami was ranked 12th best out of the 500 US cities we studied.
We proceeded to determine the number of additional supermarkets required for each of the 500 cities so that the EDE of the distance to nearest supermarket was less than or equal to:
\begin{itemize}
    \item the average EDE from all 500 US cities (2.29 km),
    \item 15 minutes (1200 m),
    \item 10 minutes (800 m),
    \item 5 minutes (400 m).
\end{itemize}
For the largest 20 cities in the US, the number of required supermarkets are shown in \autoref{tab:cities_add}.
\begin{table}[ht]
\setlength{\tabcolsep}{7pt}
\centering
    \begin{tabular}{lccccc}
    \toprule
        && \multicolumn{4}{c}{\textbf{Additional Supermarkets}} \\
        \textbf{City, State} & \textbf{Rank} & {Average} & {1200 m} & {800 m} & {400 m}   \\
        \midrule
        New York, NY      & 3   & 0  & 0   & 2   & 440  \\
        San Francisco, CA & 7   & 0  & 0   & 9   & 127  \\
        Philadelphia, PA  & 15  & 0  & 0   & 27  & 329  \\
        Washington, DC    & 17  & 0  & 0   & 15  & 182  \\
        Chicago, IL       & 20  & 0  & 0   & 43  & 553  \\
        \midrule
        Seattle, WA       & 26  & 0  & 0   & 27  & 311  \\
        Los Angeles, CA   & 56  & 0  & 8   & 186 & 2216 \\
        San Jose, CA      & 99  & 0  & 19  & 121 & --   \\
        Denver, CO        & 110 & 0  & 14  & 72  & --   \\
        San Diego, CA     & 143 & 0  & 51  & 201 & --   \\
        \midrule
        Houston, TX       & 209 & 0  & 112 & 469 & --   \\
        Charlotte, NC     & 232 & 0  & 126 & --  & --   \\
        Columbus, OH      & 260 & 0  & 81  & 270 & --   \\
        Dallas, TX        & 275 & 0  & 63  & 256 & --   \\
        Phoenix, AZ       & 321 & 0  & 80  & 358 & --   \\
        \midrule
        Indianapolis, IN  & 370 & 3  & 121 & 462 & --   \\
        San Antonio, TX   & 388 & 7  & 181 & 884 & --   \\
        Fort Worth, TX    & 449 & 8  & 103 & 371 & --   \\
        Jacksonville, FL  & 483 & 16 & 362 & --  & --   \\
        Austin, TX        & 489 & 5  & 102 & 445 & --  \\
            \bottomrule
    \end{tabular}
    \caption{
       Summary results for the 20 largest US cities sorted by rank. The last four columns indicate the number of additional supermarkets needed in each city to achieve a given level of access. ``Average'' indicates the average level of access across all 500 cities in the study, which is 2.29 km. The missing values indicate that the potential sites (Census Block Group centroids) provided to the model were insufficient for reaching the target level of access.
    }
    \label{tab:cities_add}
\end{table}
The number of supermarkets required for all 500 US cities are summarized in \autoref{fig:number-required}.

Of the 179 cities that are currently below the average level of walkable access, 73 cities require only one additional store to be on par with the average access for all 500 cities. Another 45 require two additional stores, and 43 more require between three and five new supermarkets.

To achieve an EDE value that represents a 15-minute walk, the cities in the study require on average 18 additional supermarkets (at least 8,640 supermarkets across the 487 cities where this is possible). Some cities are not far from this target. 23 cities require only one supermarket and 62 require between two and five additional supermarkets. 106 cities require between six and ten additional supermarkets, while 132 require 11 to 20 more.
With each 5-minute improvement in access, the number of new stores required increases substantially.

The missing values in \autoref{tab:cities_add} indicate that the optimization model was ``infeasible'' for that city / access-target pair.
This means that a solution could not be found based on the potential sites.
This is a limitation arising from our use of Census Block Group centroids as the potential locations for new supermarkets. In some neighborhoods, the areas of the Block Groups were too large to provide sufficient options to achieve the goal. Developing a feasible model in these cases would require including more potential sites for a more uniform coverage of the city.

\section{Discussion}
Cities and planning advocates are beginning to look towards accessible forms of urban planning in order to achieve positive sustainability and public health outcomes (e.g., the $x$-minute neighborhood).
However, the practical implications of achieving these goals have not been rigorously explored.
In this paper, we sought to answer the questions of how many supermarkets would be required to achieve accessible neighborhood goals and where they should be located so that this access is equally distributed.

We show that the scale of the retrofit required varies by city and the target level of access\autoref{fig:number-required}.
For instance, 8\% of the cities we studied already have an inequality-penalized average access of $<$1200 meters (approximately 15-minutes), and an additional 5\% are within one store of reaching that level. In contrast, more than 90\% of cities require more than 100 additional stores to achieve a five minute city target (including infeasible cities), and more than $30\%$ require more than 100 stores to reach a 10-minute target.
Unfortunately, the increase in number of stores required is not linear relative to the change in the access target. On average, we found that it would take more than 10 times as many stores to improve from 10- to 5-minute access over the number required to upgrade from 15- to 10-minute access (in the cases where those levels of access were found to be feasible).

Although these numbers are already large, they may be an underestimate of the number of stores required for a couple of reasons.
First, these estimates are based on an inequality-penalized \textit{average} for the distance people must travel to their nearest amenity.
This means that areas in the city will have to travel further than the average distance.
A more strict definition of the ``$x$-minute" neighborhood (for instance, if we attempted to ensure no one would have to travel more than $x$ minutes to their nearest amenity) would require significantly more supermarkets.
Secondly, we are not capturing ongoing sprawl and development.
The more our cities grow without careful planning, the more effort it will take to reach levels of access that are suitable for active transport.

From another angle, our results may overestimate the number of stores required.
In this study, we identify thresholds and optimize based on walking times.
In many cities, public transport will be used to extend the catchment area of amenities that are accessible by both walking and transit.
Including public transit (when it is faster than walking)  may reduce the number of supermarkets required in some cities. This type of modeling is possible if public transport data is available and if the temporal variability is appropriately considered.
Regardless of these potential discrepancies, this study provides at least an initial estimate for the scale of the challenge facing US cities seeking to enable active transport and transport choice.

Regardless, these results show that the number of supermarkets required to improve accessibility across US cities is substantial.
If cities want their residents to enjoy the wide array of public health and sustainability benefits arising from active transport and car-independent urban design they need to act.
These results speak to the urgency for this action and to the need for careful and effective planning of future development and amenity locations.

Optimizing the locations of these amenities is necessary to make this transition efficient.
The model provides optimal equitable facility locations, enabling the transition to not only be efficient, but just.
For example, for Miami, we present the optimal locations for an additional five (\autoref{fig:Miami}B) and twelve (\autoref{fig:miami-10min}) stores.

This kind of information can be used by local governments to incentivize supermarket development in particular areas and can be used by companies looking to site their next stores.
However, this study considers only the distance to the nearest supermarket, and does not consider the number of people who access the store (the demand side of access); therefore the feasibility of the supermarkets (e.g., if there is a minimum required customer base) is not considered.

When evaluating interventions such as adding amenities, it is important to consider how the benefits will be distributed across the population.
\autoref{fig:Miami}C shows how using the inequality-minimizing approach leads to notable benefits for individuals with initially poor access to supermarkets, therefore beginning to address some of the inequalities in our urban areas.
These graphs show the distribution of an intervention's benefits, and it is important that these interventions do not favor those who already have decent access.
Additionally, we observe that the gains in terms of reducing inequality did not come at the expense of decreasing the average distance (\autoref{fig:Miami}C).
Although we review the equality of the distribution, this paper does not investigate the distributional impacts between socio-demographic groups.
However, this is possible with the Kolm-Pollak EDE and is described in \cite{Logan2021-inequ}.

Although a motivation for this study is the popularity of the 15/$x$-minute concept in urban planning, these isochrone thresholds are not underlying assumptions of this method and work.
The concept of the 15, 20, or $x$-minute city nominally implies that these distances are homogeneously acceptable to residents \cite{Moreno2021-vs, Capasso_Da_Silva2019-ge, McNeil2011-zv}.
However, as is widely acknowledged, this is not the case;
one survey of walking tolerance in the Netherlands showed that only 50\% of respondents found 400 meters to be an acceptable distance to walk for food from their parked car \cite{Schaap2016-ov}.
This tolerance will likely vary significantly between cultures, climates, individuals, origins, and amenities.
To address this variability, our study evaluates a number of thresholds in order to reflect the sensitivity of the magnitude of the intervention to different accessibility targets.
Additionally, this work is not conditional on such a target.
The optimization is not based on isochrone or access-thresholds (where there are only two categories, having access or not, rather than varying levels of access). Such an approach would lead to issues with the edge-effect (i.e., residents living 15.1 minutes from their nearest supermarket considered to not have access, in contrast to their neighbour with 14.9 minutes), and may not seek to improve access further for those who are already within (or out of reach of) the access-threshold.
As we have argued in our previous work, the ultimate goal for planners should be to improve accessibility rather than to achieve some arbitrary access threshold \cite{Logan2022-mr}.

This paper presents an approach to support planners in improving their city's accessibility.
Improved accessibility (and distributional justice of this accessibility) to amenities has been directly linked to higher active transport and improved sustainability and public health outcomes \cite{Frumkin2004-yi, Pucher2010-ip, Lopez2011-cc}.
The paper is unique because it combines a new and novel measure of equality \cite{Sheriff2020-ge, Logan2021-inequ}, with an optimization model to determine where and how many amenities are required to achieve certain access goals, or simply improve accessibility.
These challenges are salient given the rise of urbanism and a wide commitment to improving distributional justice in our communities \cite{Gusdorf2007-ms, Bulkeley2014-so, Fussel2010-te}.
It is particularly salient as cities look to restore accessibility and opportunities to their communities following the pandemic and other disruptions such as natural hazards \cite{Patel2020-yq, Power2020-yb, Rose2011-gy}.
Although we apply the model to the case of supermarkets and food access, the model can be used for any amenity, from green space to healthcare to polling locations.
Ultimately, this paper and its companion methods paper \cite{horton2024scalable} lay a foundation for cities to begin to use an evidence-base to efficiently retrofit their existing development for the benefit of both the climate and their residents' health.

Ultimately, these results suggest that for many of the studied US cities, the scale of change required to achieve walkability may make retrofitting unfeasible.
This is a direct result of the urban form.
With urban sprawl continuing unabated worldwide, these US cities provide a cautionary tale from public health and sustainability perspectives.
If cities and communities are genuinely committed to enhancing their urban design to realize these benefits, the conversations must move beyond superficial commitments and focus on the role of the built environment.

\section*{Acknowledgment}
{\small This work used computing resources at the Center for Computational Mathematics, University of Colorado Denver, including the Alderaan cluster, supported by the National Science Foundation award OAC-2019089.}

\bibliography{references}

\newpage

\appendix

\include{appendix_table}

\end{document}

%% file: appendix_table.tex
\renewcommand{\arraystretch}{.8}
\begin{longtable}{lrrrrrrr}
\toprule
&&&& \multicolumn{4}{c}{Additional Supermarkets} \\
City, State & Rank & EDE (km) & Population 
& Avg & 15 min & 10 min & 5 min   \\
\midrule
\endhead
Union City, NJ                 & 1   & 579.3   & 68,186    & 0  & 0   & 0   & 3    \\
Santa Monica, CA               & 2   & 764.2   & 92,812    & 0  & 0   & 0   & 22   \\
New York, NY                   & 3   & 800.2   & 8,784,592 & 0  & 0   & 2   & 440  \\
Jersey City, NJ                & 4   & 832.7   & 291,585   & 0  & 0   & 1   & 31   \\
Cambridge, MA                  & 5   & 851.2   & 117,858   & 0  & 0   & 1   & 21   \\
Inglewood, CA                  & 6   & 948.9   & 106,817   & 0  & 0   & 2   & 44   \\
San Francisco, CA              & 7   & 968.5   & 871,136   & 0  & 0   & 9   & 127  \\
Redondo Beach, CA              & 8   & 994.1   & 71,344    & 0  & 0   & 3   & 27   \\
Berkeley, CA                   & 9   & 995.2   & 123,485   & 0  & 0   & 2   & 30   \\
Hawthorne, CA                  & 10  & 1005.8  & 87,911    & 0  & 0   & 3   & 31   \\
Somerville, MA                 & 11  & 1011.7  & 80,995    & 0  & 0   & 2   & 17   \\
Miami, FL                      & 12  & 1039.6  & 441,228   & 0  & 0   & 12  & 120  \\
Burbank, CA                    & 13  & 1041.4  & 104,508   & 0  & 0   & 3   & 52   \\
South Gate, CA                 & 14  & 1062.3  & 91,627    & 0  & 0   & 4   & 34   \\
Philadelphia, PA               & 15  & 1063.9  & 1,593,147 & 0  & 0   & 27  & 329  \\
Santa Clara, CA                & 16  & 1064.8  & 126,522   & 0  & 0   & 7   & --   \\
Washington, DC                 & 17  & 1073.5  & 684,900   & 0  & 0   & 15  & 182  \\
Providence, RI                 & 18  & 1080.2  & 189,588   & 0  & 0   & 6   & 66   \\
Long Beach, CA                 & 19  & 1089.5  & 464,262   & 0  & 0   & 15  & 164  \\
Chicago, IL                    & 20  & 1093.2  & 2,733,239 & 0  & 0   & 43  & 553  \\
Hialeah, FL                    & 21  & 1094.8  & 222,413   & 0  & 0   & 7   & --   \\
Newark, NJ                     & 22  & 1098.7  & 310,849   & 0  & 0   & 5   & 79   \\
Paterson, NJ                   & 23  & 1103.9  & 159,216   & 0  & 0   & 2   & 26   \\
Lakewood, CA                   & 24  & 1109.3  & 82,198    & 0  & 0   & 7   & --   \\
Evanston, IL                   & 25  & 1111.0  & 77,617    & 0  & 0   & 4   & 29   \\
Seattle, WA                    & 26  & 1115.5  & 726,482   & 0  & 0   & 27  & 311  \\
Mount Vernon, NY               & 27  & 1116.8  & 73,645    & 0  & 0   & 2   & 18   \\
Alexandria, VA                 & 28  & 1129.8  & 157,507   & 0  & 0   & 10  & 64   \\
Passaic, NJ                    & 29  & 1139.6  & 70,297    & 0  & 0   & 2   & 12   \\
Cicero, IL                     & 30  & 1152.9  & 85,026    & 0  & 0   & 2   & 14   \\
Whittier, CA                   & 31  & 1155.0  & 85,826    & 0  & 0   & 8   & --   \\
New Rochelle, NY               & 32  & 1167.2  & 77,661    & 0  & 0   & 6   & --   \\
Pasadena, CA                   & 33  & 1171.7  & 135,215   & 0  & 0   & 10  & --   \\
Lowell, MA                     & 34  & 1174.9  & 112,487   & 0  & 0   & 8   & --   \\
Santa Ana, CA                  & 35  & 1182.1  & 309,348   & 0  & 0   & 15  & --   \\
Sunnyvale, CA                  & 36  & 1190.5  & 154,895   & 0  & 0   & 14  & --   \\
Bellflower, CA                 & 37  & 1190.9  & 79,056    & 0  & 0   & 6   & --   \\
Lawrence, MA                   & 38  & 1191.8  & 88,522    & 0  & 0   & 4   & 28   \\
Baldwin Park, CA               & 39  & 1193.9  & 71,795    & 0  & 0   & 6   & --   \\
Huntington Beach, CA           & 40  & 1196.2  & 196,446   & 0  & 0   & 20  & --   \\
Lynwood, CA                    & 41  & 1210.2  & 67,083    & 0  & 1   & 3   & 25   \\
Daly City, CA                  & 42  & 1222.8  & 104,020   & 0  & 1   & 7   & --   \\
Yonkers, NY                    & 43  & 1224.4  & 210,151   & 0  & 1   & 13  & 105  \\
Santa Barbara, CA              & 44  & 1234.1  & 83,661    & 0  & 1   & 10  & --   \\
Pawtucket, RI                  & 45  & 1234.5  & 74,767    & 0  & 1   & 5   & 36   \\
Hartford, CT                   & 46  & 1241.3  & 119,089   & 0  & 1   & 8   & --   \\
Elizabeth, NJ                  & 47  & 1248.1  & 136,788   & 0  & 1   & 5   & 28   \\
Syracuse, NY                   & 48  & 1248.3  & 145,796   & 0  & 1   & 13  & --   \\
San Mateo, CA                  & 49  & 1249.1  & 104,267   & 0  & 1   & 8   & --   \\
Torrance, CA                   & 50  & 1252.3  & 146,799   & 0  & 1   & 14  & --   \\
Buena Park, CA                 & 51  & 1258.8  & 83,271    & 0  & 1   & 9   & --   \\
Downey, CA                     & 52  & 1259.8  & 114,624   & 0  & 1   & 12  & --   \\
El Cajon, CA                   & 53  & 1298.0  & 104,535   & 0  & 1   & 14  & --   \\
Buffalo, NY                    & 54  & 1303.8  & 274,210   & 0  & 2   & 17  & 152  \\
Erie, PA                       & 55  & 1306.6  & 92,128    & 0  & 1   & 10  & --   \\
Los Angeles, CA                & 56  & 1308.9  & 3,849,235 & 0  & 8   & 186 & 2216 \\
Boston, MA                     & 57  & 1309.6  & 670,755   & 0  & 2   & 21  & 161  \\
Fullerton, CA                  & 58  & 1314.6  & 141,320   & 0  & 2   & 22  & --   \\
Glendale, CA                   & 59  & 1315.8  & 192,354   & 0  & 2   & 12  & --   \\
Tempe, AZ                      & 60  & 1321.4  & 176,910   & 0  & 2   & 27  & --   \\
Costa Mesa, CA                 & 61  & 1336.9  & 109,236   & 0  & 2   & 12  & --   \\
Alhambra, CA                   & 62  & 1337.4  & 82,703    & 0  & 1   & 5   & --   \\
Fort Lauderdale, FL            & 63  & 1357.1  & 181,895   & 0  & 3   & 28  & --   \\
St. Louis, MO                  & 64  & 1372.9  & 298,399   & 0  & 4   & 28  & --   \\
Rochester, NY                  & 65  & 1373.9  & 208,334   & 0  & 3   & 19  & --   \\
Trenton, NJ                    & 66  & 1382.7  & 90,014    & 0  & 1   & 5   & 31   \\
Mountain View, CA              & 67  & 1388.1  & 81,213    & 0  & 1   & 8   & --   \\
Beaverton, OR                  & 68  & 1391.0  & 95,542    & 0  & 3   & 17  & --   \\
Garden Grove, CA               & 69  & 1391.3  & 171,455   & 0  & 3   & 20  & --   \\
Oakland, CA                    & 70  & 1396.6  & 432,343   & 0  & 4   & 28  & 212  \\
Alameda, CA                    & 71  & 1409.1  & 76,367    & 0  & 2   & 6   & 32   \\
Yakima, WA                     & 72  & 1419.4  & 87,551    & 0  & 3   & 19  & --   \\
Baltimore, MD                  & 73  & 1421.2  & 577,766   & 0  & 6   & 40  & 270  \\
Allentown, PA                  & 74  & 1424.1  & 123,746   & 0  & 1   & 8   & 54   \\
Ontario, CA                    & 75  & 1436.9  & 172,041   & 0  & 4   & 25  & --   \\
Bellingham, WA                 & 76  & 1439.1  & 80,221    & 0  & 3   & 18  & --   \\
Miami Beach, FL                & 77  & 1442.3  & 80,471    & 0  & 1   & 5   & 26   \\
Portland, OR                   & 78  & 1454.9  & 634,209   & 0  & 9   & 66  & --   \\
Racine, WI                     & 79  & 1455.9  & 76,989    & 0  & 2   & 9   & 69   \\
Reading, PA                    & 80  & 1462.7  & 93,921    & 0  & 1   & 5   & 25   \\
Pomona, CA                     & 81  & 1463.9  & 149,957   & 0  & 4   & 19  & --   \\
Wilmington, DE                 & 82  & 1465.8  & 70,141    & 0  & 1   & 4   & 29   \\
Everett, WA                    & 83  & 1471.1  & 103,436   & 0  & 3   & 16  & --   \\
Albany, NY                     & 84  & 1472.5  & 95,045    & 0  & 2   & 7   & --   \\
Hollywood, FL                  & 85  & 1475.1  & 152,077   & 0  & 5   & 25  & --   \\
Ventura, CA                    & 86  & 1478.8  & 107,880   & 0  & 3   & 21  & --   \\
Anaheim, CA                    & 87  & 1480.7  & 343,296   & 0  & 8   & 44  & --   \\
Norwalk, CA                    & 88  & 1503.0  & 102,487   & 0  & 2   & 9   & --   \\
Bridgeport, CT                 & 89  & 1514.8  & 147,033   & 0  & 4   & 12  & 84   \\
Brockton, MA                   & 90  & 1517.8  & 100,074   & 0  & 4   & 18  & --   \\
Arlington Heights, IL          & 91  & 1522.8  & 75,066    & 0  & 5   & 17  & --   \\
Tustin, CA                     & 92  & 1525.6  & 79,005    & 0  & 3   & 12  & --   \\
Boulder, CO                    & 93  & 1531.5  & 105,128   & 0  & 3   & 15  & --   \\
Santa Maria, CA                & 94  & 1533.3  & 108,035   & 0  & 3   & 10  & --   \\
Cleveland, OH                  & 95  & 1539.1  & 363,467   & 0  & 8   & 44  & --   \\
Compton, CA                    & 96  & 1539.8  & 95,422    & 0  & 1   & 5   & 45   \\
Turlock, CA                    & 97  & 1542.2  & 71,557    & 0  & 2   & 12  & --   \\
Upland, CA                     & 98  & 1544.4  & 77,572    & 0  & 3   & 17  & --   \\
San Jose, CA                   & 99  & 1555.4  & 993,779   & 0  & 19  & 121 & --   \\
Milwaukee, WI                  & 100 & 1569.0  & 564,921   & 0  & 13  & 56  & 438  \\
Glendale, AZ                   & 101 & 1569.5  & 241,320   & 0  & 10  & 50  & --   \\
Corona, CA                     & 102 & 1574.8  & 151,544   & 0  & 6   & 32  & --   \\
New Bedford, MA                & 103 & 1575.1  & 97,844    & 0  & 1   & 6   & 72   \\
Norfolk, VA                    & 104 & 1576.9  & 231,618   & 0  & 8   & 44  & --   \\
Spokane, WA                    & 105 & 1582.3  & 213,797   & 0  & 10  & 42  & --   \\
Mission Viejo, CA              & 106 & 1582.5  & 92,490    & 0  & 7   & 27  & --   \\
Redwood City, CA               & 107 & 1586.9  & 82,797    & 0  & 2   & 8   & --   \\
Milpitas, CA                   & 108 & 1587.9  & 79,444    & 0  & 2   & 14  & --   \\
Oxnard, CA                     & 109 & 1600.5  & 200,281   & 0  & 3   & 22  & --   \\
Denver, CO                     & 110 & 1607.5  & 705,515   & 0  & 14  & 72  & --   \\
San Leandro, CA                & 111 & 1608.5  & 90,139    & 0  & 2   & 8   & --   \\
Portland, ME                   & 112 & 1613.2  & 63,481    & 0  & 3   & 11  & --   \\
Honolulu, HI                   & 113 & 1616.4  & 341,854   & 0  & 6   & 26  & 148  \\
Bellevue, WA                   & 114 & 1629.8  & 136,046   & 0  & 9   & 36  & --   \\
Vista, CA                      & 115 & 1632.6  & 96,566    & 0  & 7   & 35  & --   \\
Largo, FL                      & 116 & 1633.8  & 78,838    & 0  & 6   & 32  & --   \\
Iowa City, IA                  & 117 & 1644.8  & 67,677    & 0  & 6   & 24  & --   \\
Miami Gardens, FL              & 118 & 1647.8  & 110,828   & 0  & 6   & 21  & --   \\
Vancouver, WA                  & 119 & 1649.4  & 172,975   & 0  & 10  & 42  & --   \\
Salinas, CA                    & 120 & 1659.3  & 161,947   & 0  & 4   & 16  & --   \\
Irvine, CA                     & 121 & 1659.7  & 302,364   & 0  & 9   & 45  & --   \\
Centennial, CO                 & 122 & 1660.1  & 100,783   & 0  & 8   & 39  & --   \\
Allen, TX                      & 123 & 1661.3  & 101,915   & 0  & 7   & 31  & --   \\
Orem, UT                       & 124 & 1665.5  & 94,784    & 0  & 5   & 19  & --   \\
Richardson, TX                 & 125 & 1674.4  & 116,636   & 0  & 6   & 30  & --   \\
Riverside, CA                  & 126 & 1675.8  & 299,063   & 0  & 13  & 65  & --   \\
Union City, CA                 & 127 & 1678.2  & 68,266    & 0  & 3   & 11  & --   \\
Gresham, OR                    & 128 & 1679.3  & 107,907   & 0  & 6   & 29  & --   \\
Newport Beach, CA              & 129 & 1685.4  & 81,967    & 0  & 8   & 24  & --   \\
Tacoma, WA                     & 130 & 1686.8  & 208,299   & 0  & 7   & 32  & --   \\
Dearborn, MI                   & 131 & 1689.7  & 107,471   & 0  & 4   & 13  & --   \\
Sacramento, CA                 & 132 & 1709.0  & 517,871   & 0  & 20  & 88  & --   \\
New Britain, CT                & 133 & 1709.2  & 71,157    & 0  & 4   & 13  & --   \\
Schenectady, NY                & 134 & 1712.9  & 65,747    & 0  & 2   & 9   & --   \\
Mesa, AZ                       & 135 & 1713.1  & 484,305   & 0  & 29  & 146 & --   \\
Pleasanton, CA                 & 136 & 1714.3  & 71,513    & 0  & 5   & 27  & --   \\
Davenport, IA                  & 137 & 1717.3  & 89,289    & 0  & 8   & 31  & --   \\
Grand Rapids, MI               & 138 & 1728.3  & 190,044   & 0  & 9   & 40  & --   \\
Pembroke Pines, FL             & 139 & 1732.7  & 157,334   & 1  & 16  & 52  & --   \\
Orange, CA                     & 140 & 1733.5  & 135,952   & 0  & 7   & 28  & --   \\
Pompano Beach, FL              & 141 & 1740.4  & 110,941   & 0  & 8   & 23  & --   \\
Bethlehem, PA                  & 142 & 1743.4  & 71,984    & 0  & 4   & 12  & --   \\
San Diego, CA                  & 143 & 1745.3  & 1,347,374 & 0  & 51  & 201 & --   \\
Hemet, CA                      & 144 & 1747.9  & 86,581    & 0  & 8   & 36  & --   \\
Carson, CA                     & 145 & 1748.7  & 94,353    & 0  & 3   & 13  & --   \\
Boise City, ID                 & 146 & 1751.8  & 214,257   & 0  & 21  & 88  & --   \\
Hammond, IN                    & 147 & 1754.1  & 74,652    & 0  & 5   & 14  & --   \\
Des Moines, IA                 & 148 & 1754.8  & 197,338   & 0  & 13  & 55  & --   \\
Fontana, CA                    & 149 & 1757.5  & 202,227   & 0  & 11  & 46  & --   \\
Kenner, LA                     & 150 & 1758.1  & 65,724    & 0  & 4   & 13  & --   \\
Boca Raton, FL                 & 151 & 1758.5  & 96,111    & 0  & 11  & 36  & --   \\
Antioch, CA                    & 152 & 1765.7  & 109,165   & 0  & 7   & 31  & --   \\
Ann Arbor, MI                  & 153 & 1767.9  & 115,073   & 0  & 6   & 21  & --   \\
Westminster, CA                & 154 & 1769.2  & 90,842    & 0  & 3   & 14  & --   \\
Fremont, CA                    & 155 & 1778.7  & 223,694   & 0  & 10  & 38  & --   \\
San Bernardino, CA             & 156 & 1779.1  & 213,770   & 0  & 9   & 42  & --   \\
Redlands, CA                   & 157 & 1780.9  & 67,242    & 0  & 4   & 19  & --   \\
Thousand Oaks, CA              & 158 & 1781.8  & 109,670   & 0  & 16  & --  & --   \\
Cary, NC                       & 159 & 1784.8  & 151,115   & 0  & 27  & --  & --   \\
Kent, WA                       & 160 & 1789.6  & 98,383    & 0  & 9   & 36  & --   \\
Plano, TX                      & 161 & 1789.9  & 278,185   & 0  & 18  & 83  & --   \\
Hayward, CA                    & 162 & 1794.8  & 159,235   & 0  & 5   & 23  & --   \\
Clearwater, FL                 & 163 & 1797.1  & 108,566   & 0  & 8   & 34  & --   \\
Rancho Cucamonga, CA           & 164 & 1800.6  & 169,450   & 0  & 12  & 52  & --   \\
Federal Way, WA                & 165 & 1809.0  & 95,059    & 0  & 6   & 28  & --   \\
Lynchburg, VA                  & 166 & 1811.9  & 60,740    & 0  & 9   & 32  & --   \\
Pasadena, TX                   & 167 & 1813.5  & 146,356   & 0  & 11  & 41  & --   \\
Deerfield Beach, FL            & 168 & 1819.1  & 86,410    & 0  & 7   & 24  & --   \\
Chandler, AZ                   & 169 & 1820.5  & 267,830   & 0  & 19  & 82  & --   \\
Hillsboro, OR                  & 170 & 1821.6  & 102,863   & 0  & 5   & 21  & --   \\
Chula Vista, CA                & 171 & 1823.7  & 267,847   & 0  & 11  & 43  & --   \\
Roanoke, VA                    & 172 & 1824.7  & 88,934    & 0  & 11  & 45  & --   \\
Cranston, RI                   & 173 & 1824.7  & 76,714    & 0  & 6   & 23  & --   \\
Eugene, OR                     & 174 & 1827.1  & 167,768   & 0  & 11  & 45  & --   \\
West Covina, CA                & 175 & 1831.4  & 107,706   & 0  & 8   & 27  & --   \\
Lake Forest, CA                & 176 & 1834.1  & 85,226    & 0  & 4   & 23  & --   \\
Worcester, MA                  & 177 & 1835.8  & 194,968   & 0  & 9   & 34  & --   \\
Renton, WA                     & 178 & 1836.3  & 101,133   & 0  & 9   & 30  & --   \\
Plantation, FL                 & 179 & 1837.9  & 90,340    & 0  & 11  & 39  & --   \\
Concord, CA                    & 180 & 1839.2  & 123,107   & 0  & 7   & 27  & --   \\
Westminster, CO                & 181 & 1839.9  & 111,754   & 0  & 10  & 40  & --   \\
Asheville, NC                  & 182 & 1845.0  & 82,804    & 0  & 18  & --  & --   \\
Pittsburgh, PA                 & 183 & 1848.2  & 295,407   & 0  & 13  & 42  & --   \\
Arlington, TX                  & 184 & 1850.3  & 381,169   & 0  & 27  & 108 & --   \\
Kennewick, WA                  & 185 & 1852.7  & 75,061    & 0  & 9   & --  & --   \\
Napa, CA                       & 186 & 1857.3  & 77,597    & 0  & 4   & 20  & --   \\
Bloomington, IN                & 187 & 1872.5  & 74,634    & 0  & 8   & 28  & --   \\
Sandy Springs, GA              & 188 & 1873.7  & 96,024    & 0  & 13  & 47  & --   \\
Gainesville, FL                & 189 & 1873.8  & 131,561   & 0  & 12  & 44  & --   \\
Medford, OR                    & 190 & 1873.8  & 80,861    & 0  & 8   & --  & --   \\
Minneapolis, MN                & 191 & 1876.3  & 426,006   & 0  & 9   & 34  & 203  \\
Salt Lake City, UT             & 192 & 1876.9  & 192,419   & 0  & 8   & 28  & --   \\
Hampton, VA                    & 193 & 1879.3  & 129,399   & 0  & 18  & --  & --   \\
Lewisville, TX                 & 194 & 1884.9  & 108,170   & 0  & 9   & 41  & --   \\
Rialto, CA                     & 195 & 1886.1  & 101,915   & 0  & 4   & 17  & --   \\
Westland, MI                   & 196 & 1889.2  & 80,668    & 0  & 7   & 28  & --   \\
Salem, OR                      & 197 & 1895.6  & 166,113   & 0  & 17  & --  & --   \\
Lauderhill, FL                 & 198 & 1896.4  & 74,339    & 0  & 6   & 14  & --   \\
Mesquite, TX                   & 199 & 1897.0  & 143,400   & 0  & 13  & 48  & --   \\
Palatine, IL                   & 200 & 1903.9  & 65,935    & 0  & 6   & 26  & --   \\
Richmond, CA                   & 201 & 1905.1  & 112,469   & 0  & 7   & 25  & --   \\
Simi Valley, CA                & 202 & 1905.6  & 117,397   & 0  & 8   & 34  & --   \\
Evansville, IN                 & 203 & 1909.0  & 105,228   & 0  & 11  & 38  & --   \\
Fresno, CA                     & 204 & 1910.0  & 526,741   & 0  & 28  & 114 & --   \\
Ogden, UT                      & 205 & 1917.3  & 83,221    & 0  & 8   & 24  & --   \\
Visalia, CA                    & 206 & 1919.8  & 137,491   & 0  & 9   & 42  & --   \\
Fall River, MA                 & 207 & 1921.6  & 90,403    & 0  & 3   & 12  & --   \\
Manteca, CA                    & 208 & 1922.9  & 75,369    & 0  & 6   & 22  & --   \\
Houston, TX                    & 209 & 1927.0  & 2,215,641 & 0  & 112 & 469 & --   \\
Springfield, MA                & 210 & 1934.3  & 148,060   & 0  & 7   & 27  & --   \\
Lakewood, CO                   & 211 & 1937.3  & 152,267   & 0  & 12  & 44  & --   \\
San Marcos, CA                 & 212 & 1937.8  & 87,032    & 0  & 9   & --  & --   \\
Longmont, CO                   & 213 & 1939.6  & 94,633    & 0  & 7   & 25  & --   \\
Champaign, IL                  & 214 & 1940.1  & 84,999    & 0  & 7   & 21  & --   \\
Elgin, IL                      & 215 & 1944.4  & 107,031   & 0  & 10  & 70  & --   \\
Chico, CA                      & 216 & 1945.6  & 95,131    & 0  & 7   & 28  & --   \\
Carrollton, TX                 & 217 & 1945.6  & 130,380   & 0  & 11  & 35  & --   \\
Sioux City, IA                 & 218 & 1947.1  & 73,047    & 0  & 6   & 19  & --   \\
Durham, NC                     & 219 & 1951.2  & 243,313   & 0  & 37  & --  & --   \\
West Palm Beach, FL            & 220 & 1953.2  & 113,725   & 0  & 9   & 37  & --   \\
Missoula, MT                   & 221 & 1955.7  & 62,101    & 0  & 10  & 31  & --   \\
Citrus Heights, CA             & 222 & 1957.3  & 86,347    & 0  & 6   & 26  & --   \\
Canton, OH                     & 223 & 1958.1  & 68,119    & 0  & 6   & 18  & --   \\
Sparks, NV                     & 224 & 1962.4  & 100,003   & 0  & 8   & 35  & --   \\
Orlando, FL                    & 225 & 1966.2  & 292,032   & 0  & 24  & 154 & --   \\
Gastonia, NC                   & 226 & 1967.3  & 64,689    & 0  & 25  & --  & --   \\
Farmington Hills, MI           & 227 & 1982.3  & 69,330    & 0  & 12  & 35  & --   \\
Lake Charles, LA               & 228 & 1985.5  & 74,947    & 0  & 13  & 51  & --   \\
Charleston, SC                 & 229 & 1986.1  & 126,534   & 0  & 40  & --  & --   \\
St. Petersburg, FL             & 230 & 1989.7  & 253,993   & 0  & 13  & 52  & --   \\
Clifton, NJ                    & 231 & 1992.7  & 88,716    & 0  & 4   & 14  & --   \\
Charlotte, NC                  & 232 & 1994.4  & 804,437   & 0  & 126 & --  & --   \\
Raleigh, NC                    & 233 & 1998.2  & 430,197   & 0  & 56  & --  & --   \\
Bolingbrook, IL                & 234 & 2002.5  & 69,312    & 0  & 8   & 29  & --   \\
Las Vegas, NV                  & 235 & 2004.5  & 623,239   & 0  & 37  & 147 & --   \\
Detroit, MI                    & 236 & 2010.0  & 625,092   & 0  & 29  & 98  & --   \\
Oceanside, CA                  & 237 & 2016.6  & 169,037   & 0  & 19  & 54  & --   \\
Greenville, NC                 & 238 & 2017.3  & 74,560    & 0  & 17  & --  & --   \\
Davie, FL                      & 239 & 2020.5  & 98,467    & 0  & 25  & --  & --   \\
Dayton, OH                     & 240 & 2024.8  & 130,685   & 0  & 12  & 36  & --   \\
Toledo, OH                     & 241 & 2031.6  & 255,596   & 0  & 20  & 64  & --   \\
Modesto, CA                    & 242 & 2036.6  & 213,408   & 0  & 10  & 39  & --   \\
Omaha, NE                      & 243 & 2045.3  & 409,421   & 0  & 30  & 102 & --   \\
Decatur, IL                    & 244 & 2047.3  & 60,573    & 0  & 13  & 57  & --   \\
Schaumburg, IL                 & 245 & 2049.5  & 77,243    & 0  & 7   & 32  & --   \\
Kenosha, WI                    & 246 & 2052.5  & 93,324    & 0  & 7   & 24  & --   \\
El Monte, CA                   & 247 & 2054.2  & 108,897   & 0  & 4   & 12  & --   \\
Scranton, PA                   & 248 & 2055.1  & 72,026    & 0  & 4   & 13  & --   \\
Merced, CA                     & 249 & 2059.9  & 84,389    & 0  & 7   & 26  & --   \\
Portsmouth, VA                 & 250 & 2060.8  & 94,416    & 0  & 10  & 32  & --   \\
Cedar Rapids, IA               & 251 & 2061.3  & 116,963   & 0  & 15  & 80  & --   \\
Lafayette, LA                  & 252 & 2064.1  & 108,594   & 0  & 20  & 58  & --   \\
Escondido, CA                  & 253 & 2070.3  & 143,784   & 0  & 17  & --  & --   \\
Garland, TX                    & 254 & 2071.7  & 238,299   & 0  & 13  & 58  & --   \\
Wilmington, NC                 & 255 & 2073.5  & 102,584   & 0  & 19  & --  & --   \\
Santa Clarita, CA              & 256 & 2079.2  & 172,029   & 0  & 15  & 61  & --   \\
Sunrise, FL                    & 257 & 2079.7  & 96,314    & 0  & 11  & --  & --   \\
Fayetteville, NC               & 258 & 2082.4  & 179,552   & 0  & 45  & --  & --   \\
Appleton, WI                   & 259 & 2087.3  & 70,697    & 0  & 7   & 20  & --   \\
Columbus, OH                   & 260 & 2090.0  & 868,417   & 0  & 81  & 270 & --   \\
Lynn, MA                       & 261 & 2090.8  & 100,574   & 0  & 3   & 8   & --   \\
Newton, MA                     & 262 & 2095.8  & 85,382    & 0  & 6   & 19  & --   \\
Tampa, FL                      & 263 & 2096.1  & 375,087   & 0  & 28  & 101 & --   \\
Lexington, KY                  & 264 & 2108.3  & 300,774   & 0  & 30  & 140 & --   \\
Southfield, MI                 & 265 & 2113.1  & 67,813    & 0  & 10  & 37  & --   \\
Livonia, MI                    & 266 & 2114.1  & 84,614    & 0  & 14  & 61  & --   \\
Wyoming, MI                    & 267 & 2115.4  & 70,740    & 0  & 10  & 36  & --   \\
Richmond, VA                   & 268 & 2124.0  & 216,832   & 0  & 15  & 46  & --   \\
Parma, OH                      & 269 & 2125.6  & 76,263    & 0  & 6   & 21  & --   \\
Rock Hill, SC                  & 270 & 2127.8  & 62,130    & 0  & 17  & --  & --   \\
Warren, MI                     & 271 & 2131.0  & 132,813   & 0  & 13  & 39  & --   \\
Boynton Beach, FL              & 272 & 2135.1  & 79,303    & 0  & 8   & 23  & --   \\
Vacaville, CA                  & 273 & 2137.5  & 98,718    & 0  & 9   & --  & --   \\
Thornton, CO                   & 274 & 2139.1  & 136,466   & 0  & 15  & 54  & --   \\
Dallas, TX                     & 275 & 2139.4  & 1,269,024 & 0  & 63  & 256 & --   \\
Green Bay, WI                  & 276 & 2140.6  & 95,923    & 0  & 10  & 37  & --   \\
Murfreesboro, TN               & 277 & 2146.6  & 133,654   & 0  & 39  & --  & --   \\
Elk Grove, CA                  & 278 & 2148.0  & 167,695   & 0  & 11  & 38  & --   \\
Moreno Valley, CA              & 279 & 2148.9  & 199,624   & 0  & 15  & 54  & --   \\
McKinney, TX                   & 280 & 2156.2  & 181,632   & 0  & 23  & 79  & --   \\
New Haven, CT                  & 281 & 2160.7  & 131,263   & 0  & 4   & 13  & --   \\
Akron, OH                      & 282 & 2161.7  & 177,615   & 0  & 17  & 52  & --   \\
Mobile, AL                     & 283 & 2163.6  & 163,412   & 0  & 37  & --  & --   \\
Colorado Springs, CO           & 284 & 2169.1  & 441,729   & 0  & 57  & 257 & --   \\
Sugar Land, TX                 & 285 & 2176.5  & 78,571    & 0  & 18  & --  & --   \\
Brooklyn Park, MN              & 286 & 2177.4  & 78,523    & 0  & 9   & 29  & --   \\
Carlsbad, CA                   & 287 & 2184.4  & 109,493   & 0  & 15  & --  & --   \\
Overland Park, KS              & 288 & 2193.5  & 181,921   & 0  & 22  & 79  & --   \\
Johns Creek, GA                & 289 & 2194.2  & 74,060    & 0  & 24  & --  & --   \\
Missouri City, TX              & 290 & 2194.7  & 71,184    & 0  & 14  & --  & --   \\
Rockford, IL                   & 291 & 2194.8  & 132,663   & 0  & 15  & 46  & --   \\
Beaumont, TX                   & 292 & 2207.6  & 102,845   & 0  & 21  & --  & --   \\
West Valley City, UT           & 293 & 2213.1  & 137,151   & 0  & 11  & 42  & --   \\
Arvada, CO                     & 294 & 2214.2  & 113,007   & 0  & 14  & 51  & --   \\
Newport News, VA               & 295 & 2214.4  & 177,661   & 0  & 19  & 62  & --   \\
Lawrence, KS                   & 296 & 2214.9  & 88,713    & 0  & 9   & 32  & --   \\
Tracy, CA                      & 297 & 2215.7  & 89,203    & 0  & 9   & 29  & --   \\
Santa Fe, NM                   & 298 & 2220.7  & 65,422    & 0  & 15  & --  & --   \\
Chino, CA                      & 299 & 2221.6  & 88,988    & 0  & 7   & 35  & --   \\
Knoxville, TN                  & 300 & 2221.7  & 165,716   & 0  & 38  & --  & --   \\
Bend, OR                       & 301 & 2232.8  & 86,052    & 0  & 14  & --  & --   \\
South Bend, IN                 & 302 & 2233.8  & 96,510    & 0  & 10  & 32  & --   \\
Flint, MI                      & 303 & 2234.8  & 73,988    & 0  & 11  & 28  & --   \\
Irving, TX                     & 304 & 2236.2  & 249,508   & 0  & 16  & 54  & --   \\
Bryan, TX                      & 305 & 2240.5  & 72,845    & 0  & 20  & --  & --   \\
Palmdale, CA                   & 306 & 2241.1  & 163,712   & 0  & 26  & --  & --   \\
Murrieta, CA                   & 307 & 2243.9  & 103,793   & 0  & 13  & 40  & --   \\
Waterbury, CT                  & 308 & 2246.2  & 105,829   & 0  & 10  & 31  & --   \\
Loveland, CO                   & 309 & 2249.5  & 70,084    & 0  & 12  & --  & --   \\
Gilbert, AZ                    & 310 & 2252.2  & 258,535   & 0  & 30  & 118 & --   \\
Lafayette, IN                  & 311 & 2253.2  & 63,860    & 0  & 8   & 25  & --   \\
Greeley, CO                    & 312 & 2253.6  & 101,471   & 0  & 12  & 45  & --   \\
Springfield, MO                & 313 & 2257.2  & 156,097   & 0  & 26  & 88  & --   \\
Manchester, NH                 & 314 & 2261.2  & 104,461   & 0  & 10  & 35  & --   \\
Greensboro, NC                 & 315 & 2263.2  & 275,555   & 0  & 64  & --  & --   \\
Temecula, CA                   & 316 & 2267.4  & 105,794   & 0  & 15  & --  & --   \\
Quincy, MA                     & 317 & 2270.8  & 99,023    & 0  & 6   & 14  & --   \\
Naperville, IL                 & 318 & 2273.1  & 143,475   & 0  & 15  & 58  & --   \\
Fort Collins, CO               & 319 & 2277.1  & 157,803   & 0  & 17  & 65  & --   \\
Avondale, AZ                   & 320 & 2277.9  & 86,924    & 0  & 13  & 36  & --   \\
Phoenix, AZ                    & 321 & 2281.7  & 1,553,053 & 0  & 80  & 358 & --   \\
Cheyenne, WY                   & 322 & 2284.8  & 60,472    & 0  & 7   & 21  & --   \\
Chattanooga, TN                & 323 & 2291.7  & 144,809   & 1  & 40  & --  & --   \\
Roswell, GA                    & 324 & 2294.3  & 79,424    & 1  & 28  & --  & --   \\
Livermore, CA                  & 325 & 2320.1  & 82,877    & 1  & 13  & --  & --   \\
Pearland, TX                   & 326 & 2325.9  & 116,526   & 1  & 22  & --  & --   \\
Stockton, CA                   & 327 & 2328.8  & 315,003   & 1  & 16  & 60  & --   \\
North Las Vegas, NV            & 328 & 2330.7  & 258,748   & 1  & 20  & 64  & --   \\
San Ramon, CA                  & 329 & 2333.3  & 77,129    & 1  & 8   & 26  & --   \\
Santa Rosa, CA                 & 330 & 2336.2  & 158,264   & 1  & 11  & 47  & --   \\
Atlanta, GA                    & 331 & 2337.6  & 460,547   & 1  & 33  & 123 & --   \\
Peoria, IL                     & 332 & 2339.8  & 100,459   & 1  & 16  & 67  & --   \\
Indio, CA                      & 333 & 2342.7  & 82,254    & 1  & 14  & 38  & --   \\
Peoria, AZ                     & 334 & 2349.6  & 180,203   & 1  & 40  & --  & --   \\
Tyler, TX                      & 335 & 2353.4  & 91,719    & 1  & 25  & --  & --   \\
Pueblo, CO                     & 336 & 2355.1  & 105,150   & 1  & 12  & 32  & --   \\
Little Rock, AR                & 337 & 2355.7  & 173,674   & 1  & 54  & --  & --   \\
Rochester, MN                  & 338 & 2355.9  & 99,740    & 1  & 21  & --  & --   \\
Fayetteville, AR               & 339 & 2356.7  & 78,871    & 1  & 27  & --  & --   \\
Warner Robins, GA              & 340 & 2361.9  & 70,494    & 1  & 29  & --  & --   \\
Mount Pleasant, SC             & 341 & 2365.5  & 71,880    & 1  & --  & --  & --   \\
Folsom, CA                     & 342 & 2371.2  & 76,978    & 1  & 16  & --  & --   \\
Norwalk, CT                    & 343 & 2376.4  & 81,524    & 1  & 8   & 26  & --   \\
Springdale, AR                 & 344 & 2377.4  & 70,967    & 1  & 18  & --  & --   \\
Camden, NJ                     & 345 & 2388.6  & 71,190    & 1  & 3   & 7   & 36   \\
Independence, MO               & 346 & 2390.7  & 106,725   & 1  & 24  & 76  & --   \\
Meridian, ID                   & 347 & 2401.4  & 98,984    & 1  & 13  & 44  & --   \\
Scottsdale, AZ                 & 348 & 2401.8  & 205,084   & 1  & 62  & --  & --   \\
Muncie, IN                     & 349 & 2403.7  & 57,695    & 1  & 8   & 23  & --   \\
Tallahassee, FL                & 350 & 2407.8  & 181,092   & 1  & 42  & --  & --   \\
Deltona, FL                    & 351 & 2408.8  & 82,602    & 1  & 18  & --  & --   \\
Waukegan, IL                   & 352 & 2409.4  & 85,092    & 1  & 8   & 26  & --   \\
Coral Springs, FL              & 353 & 2429.5  & 133,606   & 1  & --  & --  & --   \\
Sioux Falls, SD                & 354 & 2433.5  & 172,902   & 1  & 23  & 100 & --   \\
Troy, MI                       & 355 & 2440.9  & 76,482    & 1  & 14  & --  & --   \\
Stamford, CT                   & 356 & 2446.0  & 116,053   & 1  & 8   & 20  & --   \\
Rochester Hills, MI            & 357 & 2446.7  & 63,088    & 1  & 18  & --  & --   \\
Birmingham, AL                 & 358 & 2467.0  & 180,448   & 2  & 37  & 115 & --   \\
Youngstown, OH                 & 359 & 2473.3  & 51,629    & 1  & 11  & 30  & --   \\
Melbourne, FL                  & 360 & 2487.8  & 80,608    & 1  & 22  & --  & --   \\
Cincinnati, OH                 & 361 & 2487.8  & 298,967   & 1  & 27  & 83  & --   \\
Chino Hills, CA                & 362 & 2500.0  & 67,639    & 1  & 11  & --  & --   \\
Columbia, SC                   & 363 & 2503.9  & 121,572   & 1  & 31  & --  & --   \\
Bakersfield, CA                & 364 & 2505.5  & 391,761   & 1  & 28  & 113 & --   \\
Madison, WI                    & 365 & 2508.6  & 255,295   & 1  & 19  & 64  & --   \\
Macon, GA                      & 366 & 2519.0  & 75,294    & 1  & 17  & 56  & --   \\
Sterling Heights, MI           & 367 & 2524.0  & 127,785   & 1  & 17  & 54  & --   \\
Spokane Valley, WA             & 368 & 2525.4  & 92,329    & 1  & 15  & --  & --   \\
Aurora, IL                     & 369 & 2528.4  & 171,214   & 1  & 17  & 59  & --   \\
Indianapolis, IN               & 370 & 2539.5  & 788,869   & 3  & 121 & 462 & --   \\
Fort Wayne, IN                 & 371 & 2539.8  & 237,047   & 2  & 43  & --  & --   \\
Vallejo, CA                    & 372 & 2546.7  & 121,109   & 1  & 10  & 31  & --   \\
Louisville, KY                 & 373 & 2557.1  & 567,257   & 3  & 89  & 337 & --   \\
Kalamazoo, MI                  & 374 & 2563.1  & 67,303    & 1  & 8   & 24  & --   \\
Nampa, ID                      & 375 & 2567.7  & 92,738    & 1  & 15  & --  & --   \\
Nashville, TN                  & 376 & 2585.7  & 588,971   & 5  & 156 & --  & --   \\
Roseville, CA                  & 377 & 2587.2  & 139,356   & 1  & 13  & 52  & --   \\
Frisco, TX                     & 378 & 2594.9  & 187,075   & 1  & 27  & --  & --   \\
Provo, UT                      & 379 & 2595.6  & 110,239   & 1  & 8   & 25  & --   \\
Winston-Salem, NC              & 380 & 2610.8  & 202,683   & 2  & 89  & --  & --   \\
Olathe, KS                     & 381 & 2611.5  & 127,633   & 1  & 21  & 65  & --   \\
Clovis, CA                     & 382 & 2617.8  & 110,618   & 1  & 8   & 30  & --   \\
Baytown, TX                    & 383 & 2622.8  & 78,510    & 2  & 17  & --  & --   \\
Chesapeake, VA                 & 384 & 2624.5  & 207,651   & 3  & --  & --  & --   \\
Lansing, MI                    & 385 & 2630.2  & 106,227   & 2  & 13  & 37  & --   \\
Lancaster, CA                  & 386 & 2630.8  & 163,417   & 1  & 16  & --  & --   \\
Albuquerque, NM                & 387 & 2639.0  & 549,255   & 3  & 48  & 156 & --   \\
San Antonio, TX                & 388 & 2641.9  & 1,381,080 & 7  & 181 & 884 & --   \\
Henderson, NV                  & 389 & 2644.6  & 303,432   & 2  & 41  & --  & --   \\
North Charleston, SC           & 390 & 2646.5  & 102,806   & 2  & --  & --  & --   \\
Burlington, VT                 & 391 & 2646.5  & 42,520    & 1  & 3   & 8   & --   \\
Layton, UT                     & 392 & 2650.0  & 72,740    & 2  & 12  & --  & --   \\
Charleston, WV                 & 393 & 2668.3  & 36,278    & 2  & 15  & --  & --   \\
Yuma, AZ                       & 394 & 2684.5  & 87,450    & 1  & 10  & 51  & --   \\
Memphis, TN                    & 395 & 2692.6  & 585,130   & 3  & 77  & 270 & --   \\
Lubbock, TX                    & 396 & 2702.2  & 238,590   & 3  & 28  & 92  & --   \\
Danbury, CT                    & 397 & 2702.3  & 68,197    & 1  & 13  & 45  & --   \\
League City, TX                & 398 & 2706.0  & 107,397   & 1  & 20  & 58  & --   \\
Jonesboro, AR                  & 399 & 2711.0  & 55,667    & 2  & 22  & --  & --   \\
St. Paul, MN                   & 400 & 2725.5  & 306,600   & 1  & 14  & 39  & --   \\
El Paso, TX                    & 401 & 2731.1  & 640,366   & 3  & 61  & 198 & --   \\
Corpus Christi, TX             & 402 & 2733.7  & 295,863   & 3  & 44  & 143 & --   \\
Lakeland, FL                   & 403 & 2734.2  & 102,468   & 1  & 19  & --  & --   \\
Springfield, IL                & 404 & 2736.7  & 101,948   & 2  & 28  & --  & --   \\
Waterloo, IA                   & 405 & 2747.1  & 57,317    & 2  & 10  & 31  & --   \\
Redding, CA                    & 406 & 2749.5  & 71,725    & 2  & 41  & --  & --   \\
Topeka, KS                     & 407 & 2759.0  & 115,006   & 2  & 19  & 60  & --   \\
McAllen, TX                    & 408 & 2775.0  & 131,993   & 2  & 15  & 53  & --   \\
High Point, NC                 & 409 & 2775.6  & 97,657    & 2  & 31  & --  & --   \\
Kansas City, MO                & 410 & 2785.6  & 453,839   & 3  & 64  & 209 & --   \\
Grand Prairie, TX              & 411 & 2788.0  & 188,187   & 2  & 27  & --  & --   \\
Fairfield, CA                  & 412 & 2793.0  & 112,955   & 1  & 14  & --  & --   \\
Fargo, ND                      & 413 & 2794.6  & 117,190   & 2  & 13  & 38  & --   \\
Mission, TX                    & 414 & 2804.4  & 79,832    & 2  & 16  & --  & --   \\
Longview, TX                   & 415 & 2804.6  & 67,431    & 2  & 20  & --  & --   \\
Warwick, RI                    & 416 & 2807.4  & 73,689    & 2  & 14  & 50  & --   \\
Savannah, GA                   & 417 & 2811.0  & 137,251   & 2  & 24  & --  & --   \\
Montgomery, AL                 & 418 & 2858.6  & 171,788   & 4  & 51  & --  & --   \\
Baton Rouge, LA                & 419 & 2882.1  & 210,391   & 3  & 38  & 116 & --   \\
Broken Arrow, OK               & 420 & 2885.9  & 93,647    & 3  & 31  & --  & --   \\
Sandy, UT                      & 421 & 2889.2  & 90,722    & 1  & 12  & 45  & --   \\
College Station, TX            & 422 & 2896.3  & 113,469   & 2  & 16  & --  & --   \\
Bloomington, MN                & 423 & 2922.5  & 81,089    & 1  & 14  & 37  & --   \\
Lee's Summit, MO               & 424 & 2928.4  & 87,205    & 3  & 28  & --  & --   \\
Plymouth, MN                   & 425 & 2955.2  & 64,676    & 2  & 18  & --  & --   \\
Columbus, GA                   & 426 & 2961.1  & 182,045   & 4  & --  & --  & --   \\
San Angelo, TX                 & 427 & 2969.5  & 87,305    & 1  & 13  & 46  & --   \\
Hesperia, CA                   & 428 & 2969.7  & 87,113    & 2  & 17  & --  & --   \\
Auburn, WA                     & 429 & 2970.6  & 79,338    & 2  & 11  & 31  & --   \\
Nashua, NH                     & 430 & 2976.3  & 80,343    & 1  & 12  & 37  & --   \\
Las Cruces, NM                 & 431 & 2988.1  & 103,326   & 2  & 21  & --  & --   \\
Tulsa, OK                      & 432 & 3022.7  & 380,249   & 3  & 45  & 163 & --   \\
Fort Smith, AR                 & 433 & 3036.1  & 77,013    & 2  & 18  & --  & --   \\
Round Rock, TX                 & 434 & 3041.8  & 104,897   & 3  & 20  & --  & --   \\
Hoover, AL                     & 435 & 3056.6  & 76,966    & 6  & 37  & --  & --   \\
Victorville, CA                & 436 & 3079.8  & 127,018   & 4  & 27  & --  & --   \\
Killeen, TX                    & 437 & 3098.4  & 139,105   & 2  & 19  & 77  & --   \\
Gary, IN                       & 438 & 3098.6  & 62,903    & 2  & 11  & 34  & --   \\
Fishers, IN                    & 439 & 3109.7  & 87,028    & 4  & --  & --  & --   \\
Menifee, CA                    & 440 & 3113.0  & 90,654    & 2  & 22  & --  & --   \\
Athens, GA                     & 441 & 3139.5  & 95,851    & 2  & 52  & --  & --   \\
Edmond, OK                     & 442 & 3147.7  & 72,976    & 3  & 27  & --  & --   \\
Reno, NV                       & 443 & 3157.7  & 242,706   & 3  & 33  & 106 & --   \\
Kansas City, KS                & 444 & 3167.0  & 128,020   & 4  & 26  & 81  & --   \\
Concord, NC                    & 445 & 3169.1  & 81,727    & 4  & --  & --  & --   \\
Billings, MT                   & 446 & 3199.0  & 104,701   & 2  & 25  & --  & --   \\
Wichita, KS                    & 447 & 3204.4  & 360,255   & 10 & 71  & --  & --   \\
Albany, GA                     & 448 & 3233.7  & 57,411    & 3  & 24  & --  & --   \\
Fort Worth, TX                 & 449 & 3310.1  & 865,707   & 8  & 103 & 371 & --   \\
Suffolk, VA                    & 450 & 3315.3  & 57,328    & 4  & 26  & --  & --   \\
Laredo, TX                     & 451 & 3340.9  & 242,729   & 4  & 26  & 78  & --   \\
Joliet, IL                     & 452 & 3344.8  & 141,960   & 3  & 24  & --  & --   \\
Aurora, CO                     & 453 & 3344.9  & 374,055   & 1  & 30  & 138 & --   \\
Augusta, GA                    & 454 & 3346.5  & 156,034   & 9  & 84  & --  & --   \\
Norman, OK                     & 455 & 3373.0  & 104,271   & 2  & 19  & 58  & --   \\
O'Fallon, MO                   & 456 & 3373.9  & 83,465    & 4  & 30  & --  & --   \\
Clarksville, TN                & 457 & 3376.7  & 132,241   & 7  & --  & --  & --   \\
Apple Valley, CA               & 458 & 3378.7  & 62,164    & 4  & 25  & --  & --   \\
Carmel, IN                     & 459 & 3411.7  & 75,908    & 2  & 32  & --  & --   \\
Gulfport, MS                   & 460 & 3428.7  & 59,559    & 4  & 23  & --  & --   \\
Waukesha, WI                   & 461 & 3443.5  & 65,476    & 1  & 8   & 25  & --   \\
Rapid City, SD                 & 462 & 3475.2  & 55,013    & 3  & 20  & --  & --   \\
Oklahoma City, OK              & 463 & 3486.2  & 585,298   & 12 & 114 & 427 & --   \\
Bloomington, IL                & 464 & 3514.0  & 73,622    & 3  & 15  & --  & --   \\
Perris, CA                     & 465 & 3534.6  & 74,679    & 2  & 9   & 25  & --   \\
West Jordan, UT                & 466 & 3544.3  & 113,136   & 2  & 14  & 57  & --   \\
Jackson, MS                    & 467 & 3545.0  & 132,215   & 7  & 39  & --  & --   \\
Pharr, TX                      & 468 & 3584.1  & 77,075    & 2  & 15  & 38  & --   \\
St. Joseph, MO                 & 469 & 3615.3  & 61,275    & 2  & 14  & 38  & --   \\
Palm Bay, FL                   & 470 & 3681.5  & 106,342   & 5  & 43  & --  & --   \\
Edinburg, TX                   & 471 & 3729.2  & 79,803    & 3  & 19  & --  & --   \\
New Orleans, LA                & 472 & 3741.0  & 377,281   & 1  & 15  & 60  & --   \\
Midland, TX                    & 473 & 3741.5  & 115,533   & 2  & 16  & 73  & --   \\
Lincoln, NE                    & 474 & 3784.6  & 270,110   & 3  & 26  & 86  & --   \\
Tuscaloosa, AL                 & 475 & 3862.5  & 84,424    & 1  & 21  & --  & --   \\
Amarillo, TX                   & 476 & 3900.6  & 182,518   & 2  & 25  & 73  & --   \\
Columbia, MO                   & 477 & 3905.6  & 107,706   & 4  & 37  & --  & --   \\
Huntsville, AL                 & 478 & 3922.4  & 180,633   & 8  & --  & --  & --   \\
Brownsville, TX                & 479 & 3989.6  & 178,648   & 5  & 27  & 78  & --   \\
Wichita Falls, TX              & 480 & 4013.7  & 90,894    & 2  & 16  & 45  & --   \\
Odessa, TX                     & 481 & 4087.1  & 102,381   & 2  & 11  & 36  & --   \\
Shreveport, LA                 & 482 & 4100.1  & 163,611   & 9  & 61  & --  & --   \\
Jacksonville, FL               & 483 & 4102.4  & 834,225   & 16 & 362 & --  & --   \\
Port St. Lucie, FL             & 484 & 4118.5  & 192,631   & 10 & --  & --  & --   \\
Cape Coral, FL                 & 485 & 4204.2  & 179,013   & 9  & 88  & --  & --   \\
Surprise, AZ                   & 486 & 4249.6  & 136,583   & 3  & 28  & --  & --   \\
Virginia Beach, VA             & 487 & 4680.5  & 417,731   & 1  & 58  & --  & --   \\
St. George, UT                 & 488 & 4705.5  & 81,834    & 6  & --  & --  & --   \\
Austin, TX                     & 489 & 4846.7  & 893,947   & 5  & 102 & 445 & --   \\
Tucson, AZ                     & 490 & 4892.2  & 508,571   & 1  & 29  & 132 & --   \\
Waco, TX                       & 491 & 4913.4  & 125,775   & 7  & 69  & --  & --   \\
Palm Coast, FL                 & 492 & 5050.1  & 77,271    & 11 & --  & --  & --   \\
Jacksonville, NC               & 493 & 5280.3  & 60,367    & 2  & 24  & --  & --   \\
Duluth, MN                     & 494 & 5486.4  & 66,602    & 3  & 12  & 26  & --   \\
Rio Rancho, NM                 & 495 & 5950.8  & 92,607    & 6  & 36  & --  & --   \\
Lawton, OK                     & 496 & 6334.7  & 80,836    & 4  & 30  & --  & --   \\
Abilene, TX                    & 497 & 7097.8  & 109,596   & 5  & 30  & --  & --   \\
Denton, TX                     & 498 & 7297.7  & 121,649   & 3  & 22  & --  & --   \\
Miramar, FL                    & 499 & 11850.7 & 132,450   & 16 & --  & --  & --   \\
Anchorage, AK                  & 500 & 34723.5 & 193,910   & 2  & 29  & 114 & --   \\
\end{longtable}

%% file: supermarketoptimization.bbl
\begin{thebibliography}{10}
\expandafter\ifx\csname url\endcsname\relax
  \def\url#1{\texttt{#1}}\fi
\expandafter\ifx\csname urlprefix\endcsname\relax\def\urlprefix{URL }\fi
\providecommand{\bibinfo}[2]{#2}
\providecommand{\eprint}[2][]{\url{#2}}

\bibitem{McCormack2008-nl}
\bibinfo{author}{McCormack, G.~R.}, \bibinfo{author}{Giles-Corti, B.} \& \bibinfo{author}{Bulsara, M.}
\newblock \bibinfo{title}{The relationship between destination proximity, destination mix and physical activity behaviors}.
\newblock \emph{\bibinfo{journal}{Prev. Med.}} \textbf{\bibinfo{volume}{46}}, \bibinfo{pages}{33--40} (\bibinfo{year}{2008}).

\bibitem{Celis-Morales2017-mn}
\bibinfo{author}{Celis-Morales, C.~A.} \emph{et~al.}
\newblock \bibinfo{title}{Association between active commuting and incident cardiovascular disease, cancer, and mortality: prospective cohort study}.
\newblock \emph{\bibinfo{journal}{BMJ}} \textbf{\bibinfo{volume}{357}}, \bibinfo{pages}{j1456} (\bibinfo{year}{2017}).

\bibitem{Lindsay2011-ys}
\bibinfo{author}{Lindsay, G.}, \bibinfo{author}{Macmillan, A.} \& \bibinfo{author}{Woodward, A.}
\newblock \bibinfo{title}{Moving urban trips from cars to bicycles: impact on health and emissions}.
\newblock \emph{\bibinfo{journal}{Aust. N. Z. J. Public Health}} \textbf{\bibinfo{volume}{35}}, \bibinfo{pages}{54--60} (\bibinfo{year}{2011}).

\bibitem{Stevenson2016-yb}
\bibinfo{author}{Stevenson, M.} \emph{et~al.}
\newblock \bibinfo{title}{Land use, transport, and population health: estimating the health benefits of compact cities}.
\newblock \emph{\bibinfo{journal}{Lancet}} \textbf{\bibinfo{volume}{388}}, \bibinfo{pages}{2925--2935} (\bibinfo{year}{2016}).

\bibitem{Mueller2020-bl}
\bibinfo{author}{Mueller, N.} \emph{et~al.}
\newblock \bibinfo{title}{Changing the urban design of cities for health: The superblock model}.
\newblock \emph{\bibinfo{journal}{Environ. Int.}} \textbf{\bibinfo{volume}{134}}, \bibinfo{pages}{105132} (\bibinfo{year}{2020}).

\bibitem{Staricco2022-zk}
\bibinfo{author}{Staricco, L.} \& \bibinfo{author}{Brovarone, E.~V.}
\newblock \bibinfo{title}{Livable neighborhoods for sustainable cities: Insights from barcelona}.
\newblock \emph{\bibinfo{journal}{Transportation Research Procedia}} \textbf{\bibinfo{volume}{60}}, \bibinfo{pages}{354--361} (\bibinfo{year}{2022}).

\bibitem{Leyden2003-dg}
\bibinfo{author}{Leyden, K.~M.}
\newblock \bibinfo{title}{Social capital and the built environment: the importance of walkable neighborhoods}.
\newblock \emph{\bibinfo{journal}{Am. J. Public Health}} \textbf{\bibinfo{volume}{93}}, \bibinfo{pages}{1546--1551} (\bibinfo{year}{2003}).

\bibitem{Lopez2011-cc}
\bibinfo{author}{Lopez, R.~P.}
\newblock \emph{\bibinfo{title}{The Built Environment and Public Health}} (\bibinfo{publisher}{John Wiley \& Sons}, \bibinfo{year}{2011}).

\bibitem{Brand2021-qw}
\bibinfo{author}{Brand, C.} \emph{et~al.}
\newblock \bibinfo{title}{The climate change mitigation impacts of active travel: Evidence from a longitudinal panel study in seven european cities}.
\newblock \emph{\bibinfo{journal}{Glob. Environ. Change}} \textbf{\bibinfo{volume}{67}}, \bibinfo{pages}{102224} (\bibinfo{year}{2021}).

\bibitem{Iea2020-is}
\bibinfo{author}{{IEA}}.
\newblock \bibinfo{title}{Tracking transport 2020}.
\newblock \bibinfo{howpublished}{\url{https://www.iea.org/reports/tracking-transport-2020}} (\bibinfo{year}{2020}).
\newblock \bibinfo{note}{Accessed: 2021-12-18}.

\bibitem{Dempsey2011-og}
\bibinfo{author}{Dempsey, N.}, \bibinfo{author}{Bramley, G.}, \bibinfo{author}{Power, S.} \& \bibinfo{author}{Brown, C.}
\newblock \bibinfo{title}{The social dimension of sustainable development: Defining urban social sustainability}.
\newblock \emph{\bibinfo{journal}{Sust. Dev.}} \textbf{\bibinfo{volume}{19}}, \bibinfo{pages}{289--300} (\bibinfo{year}{2011}).

\bibitem{Jacobs1961-po}
\bibinfo{author}{Jacobs, J.}
\newblock \emph{\bibinfo{title}{The Death and Life of Great American Cities}} (\bibinfo{publisher}{Vintage Books}, \bibinfo{address}{New York}, \bibinfo{year}{1961}).

\bibitem{Hao2022-tp}
\bibinfo{author}{Hao, H.} \& \bibinfo{author}{Wang, Y.}
\newblock \bibinfo{title}{Disentangling relations between urban form and urban accessibility for resilience to extreme weather and climate events}.
\newblock \emph{\bibinfo{journal}{Landsc. Urban Plan.}} \textbf{\bibinfo{volume}{220}}, \bibinfo{pages}{104352} (\bibinfo{year}{2022}).

\bibitem{Logan2021-inequ}
\bibinfo{author}{Logan, T.~M.}, \bibinfo{author}{Anderson, M.~J.}, \bibinfo{author}{Williams, T.~G.} \& \bibinfo{author}{Conrow, L.}
\newblock \bibinfo{title}{Measuring inequalities in urban systems: An approach for evaluating the distribution of amenities and burdens}.
\newblock \emph{\bibinfo{journal}{Comput. Environ. Urban Syst.}} \textbf{\bibinfo{volume}{86}}, \bibinfo{pages}{101590} (\bibinfo{year}{2021}).

\bibitem{Anderson2022-vr}
\bibinfo{author}{Anderson, M.~J.}, \bibinfo{author}{Kiddle, D. A.~F.} \& \bibinfo{author}{Logan, T.~M.}
\newblock \bibinfo{title}{The underestimated role of the transportation network: Improving disaster \& community resilience}.
\newblock \emph{\bibinfo{journal}{Transp. Res. Part D: Trans. Environ.}} \textbf{\bibinfo{volume}{106}}, \bibinfo{pages}{103218} (\bibinfo{year}{2022}).

\bibitem{Sajjad2021-vq}
\bibinfo{author}{Sajjad, M.}, \bibinfo{author}{Chan, J. C.~L.} \& \bibinfo{author}{Chopra, S.~S.}
\newblock \bibinfo{title}{Rethinking disaster resilience in high-density cities: Towards an urban resilience knowledge system}.
\newblock \emph{\bibinfo{journal}{Sustainable Cities and Society}} \textbf{\bibinfo{volume}{69}}, \bibinfo{pages}{102850} (\bibinfo{year}{2021}).

\bibitem{Logan2022-mr}
\bibinfo{author}{Logan, T.~M.} \emph{et~al.}
\newblock \bibinfo{title}{The x-minute city: Measuring the 10, 15, 20-minute city and an evaluation of its use for sustainable urban design}.
\newblock \emph{\bibinfo{journal}{Cities}} \textbf{\bibinfo{volume}{131}}, \bibinfo{pages}{103924} (\bibinfo{year}{2022}).

\bibitem{Wu2021-nx}
\bibinfo{author}{Wu, H.} \emph{et~al.}
\newblock \bibinfo{title}{Urban access across the globe: an international comparison of different transport modes}.
\newblock \emph{\bibinfo{journal}{npj Urban Sustainability}} \textbf{\bibinfo{volume}{1}}, \bibinfo{pages}{1--9} (\bibinfo{year}{2021}).

\bibitem{Pucher2010-ip}
\bibinfo{author}{Pucher, J.} \& \bibinfo{author}{Buehler, R.}
\newblock \bibinfo{title}{Walking and cycling for healthy cities}.
\newblock \emph{\bibinfo{journal}{Built Environ.}} \textbf{\bibinfo{volume}{36}}, \bibinfo{pages}{391--414} (\bibinfo{year}{2010}).

\bibitem{Behnisch2022-xa}
\bibinfo{author}{Behnisch, M.}, \bibinfo{author}{Kr{\"u}ger, T.} \& \bibinfo{author}{Jaeger, J. A.~G.}
\newblock \bibinfo{title}{Rapid rise in urban sprawl: Global hotspots and trends since 1990}.
\newblock \emph{\bibinfo{journal}{PLOS Sustainability and Transformation}} \textbf{\bibinfo{volume}{1}}, \bibinfo{pages}{e0000034} (\bibinfo{year}{2022}).

\bibitem{C40_Cities2020-bm}
\bibinfo{author}{{C40 Cities}}.
\newblock \bibinfo{title}{{C40} mayors' agenda for a green and just recovery}.
\newblock \bibinfo{type}{Tech. Rep.} (\bibinfo{year}{2020}).

\bibitem{Marino2012-lq}
\bibinfo{author}{Marino, E.} \& \bibinfo{author}{Ribot, J.}
\newblock \bibinfo{title}{Special issue introduction: Adding insult to injury: Climate change and the inequities of climate intervention}.
\newblock \emph{\bibinfo{journal}{Glob. Environ. Change}} \textbf{\bibinfo{volume}{22}}, \bibinfo{pages}{323--328} (\bibinfo{year}{2012}).

\bibitem{Wilson2008-yk}
\bibinfo{author}{Wilson, S.}, \bibinfo{author}{Hutson, M.} \& \bibinfo{author}{Mujahid, M.}
\newblock \bibinfo{title}{How planning and zoning contribute to inequitable development, neighborhood health, and environmental injustice}.
\newblock \emph{\bibinfo{journal}{Environ. Justice}} \textbf{\bibinfo{volume}{1}}, \bibinfo{pages}{211--216} (\bibinfo{year}{2008}).

\bibitem{Calvin2017-ja}
\bibinfo{author}{Calvin, K.} \emph{et~al.}
\newblock \bibinfo{title}{The {SSP4}: A world of deepening inequality}.
\newblock \emph{\bibinfo{journal}{Glob. Environ. Change}} \textbf{\bibinfo{volume}{42}}, \bibinfo{pages}{284--296} (\bibinfo{year}{2017}).

\bibitem{Fussel2010-te}
\bibinfo{author}{F{\"u}ssel, H.-M.}
\newblock \bibinfo{title}{How inequitable is the global distribution of responsibility, capability, and vulnerability to climate change: A comprehensive indicator-based assessment}.
\newblock \emph{\bibinfo{journal}{Glob. Environ. Change}} \textbf{\bibinfo{volume}{20}}, \bibinfo{pages}{597--611} (\bibinfo{year}{2010}).

\bibitem{Bulkeley2014-so}
\bibinfo{author}{Bulkeley, H.}, \bibinfo{author}{Edwards, G. A.~S.} \& \bibinfo{author}{Fuller, S.}
\newblock \bibinfo{title}{Contesting climate justice in the city: Examining politics and practice in urban climate change experiments}.
\newblock \emph{\bibinfo{journal}{Glob. Environ. Change}} \textbf{\bibinfo{volume}{25}}, \bibinfo{pages}{31--40} (\bibinfo{year}{2014}).

\bibitem{White2021-jd}
\bibinfo{author}{White, A.~G.}, \bibinfo{author}{Guikema, S.~D.} \& \bibinfo{author}{Logan, T.~M.}
\newblock \bibinfo{title}{Urban population characteristics and their correlation with historic discriminatory housing practices}.
\newblock \emph{\bibinfo{journal}{Appl. Geogr.}} \textbf{\bibinfo{volume}{132}}, \bibinfo{pages}{102445} (\bibinfo{year}{2021}).

\bibitem{Low2013-yx}
\bibinfo{author}{Low, S.}
\newblock \bibinfo{title}{Public space and diversity: Distributive, procedural and interactional justice for parks}.
\newblock \emph{\bibinfo{journal}{The Ashgate research companion to planning and culture. Ashgate Publishing, Surrey}} \bibinfo{pages}{295--310} (\bibinfo{year}{2013}).

\bibitem{Sheriff2020-ge}
\bibinfo{author}{Sheriff, G.} \& \bibinfo{author}{Maguire, K.}
\newblock \bibinfo{title}{Health risk, inequality indexes, and environmental justice}.
\newblock \emph{\bibinfo{journal}{Risk Anal.}}  (\bibinfo{year}{2020}).

\bibitem{horton2024scalable}
\bibinfo{author}{Horton, D.}, \bibinfo{author}{Logan, T.}, \bibinfo{author}{Murrell, J.}, \bibinfo{author}{Skipper, D.} \& \bibinfo{author}{Speakman, E.}
\newblock \bibinfo{title}{A scalable approach to equitable facility location} (\bibinfo{year}{2024}).
\newblock \eprint{2401.15452}.

\bibitem{foodinescurity}
\bibinfo{title}{{Economic Research Service U.S. DEPARTMENT OF AGRICULTURE} definitions of food security}.
\newblock \bibinfo{howpublished}{\url{https://www.ers.usda.gov/topics/food-nutrition-assistance/food-security-in-the-us/definitions-of-food-security.aspx}}.
\newblock \bibinfo{note}{Accessed: 2021-12-15}.

\bibitem{Kolak2018-az}
\bibinfo{author}{Kolak, M.} \emph{et~al.}
\newblock \bibinfo{title}{Urban foodscape trends: Disparities in healthy food access in chicago, 2007-2014}.
\newblock \emph{\bibinfo{journal}{Health Place}} \textbf{\bibinfo{volume}{52}}, \bibinfo{pages}{231--239} (\bibinfo{year}{2018}).

\bibitem{Garcia2020-xt}
\bibinfo{author}{Garcia, X.}, \bibinfo{author}{Garcia-Sierra, M.} \& \bibinfo{author}{Domene, E.}
\newblock \bibinfo{title}{Spatial inequality and its relationship with local food environments: The case of barcelona}.
\newblock \emph{\bibinfo{journal}{Appl. Geogr.}} \textbf{\bibinfo{volume}{115}}, \bibinfo{pages}{102140} (\bibinfo{year}{2020}).

\bibitem{Apparicio2007-di}
\bibinfo{author}{Apparicio, P.}, \bibinfo{author}{Cloutier, M.-S.} \& \bibinfo{author}{Shearmur, R.}
\newblock \bibinfo{title}{The case of montr{\'e}al's missing food deserts: Evaluation of accessibility to food supermarkets}.
\newblock \emph{\bibinfo{journal}{Int. J. Health Geogr.}} \textbf{\bibinfo{volume}{6}}, \bibinfo{pages}{4} (\bibinfo{year}{2007}).

\bibitem{Walker2010-ch}
\bibinfo{author}{Walker, R.~E.}, \bibinfo{author}{Keane, C.~R.} \& \bibinfo{author}{Burke, J.~G.}
\newblock \bibinfo{title}{Disparities and access to healthy food in the united states: A review of food deserts literature}.
\newblock \emph{\bibinfo{journal}{Health Place}} \textbf{\bibinfo{volume}{16}}, \bibinfo{pages}{876--884} (\bibinfo{year}{2010}).

\bibitem{Krenichyn2006-ve}
\bibinfo{author}{Krenichyn, K.}
\newblock \bibinfo{title}{`the only place to go and be in the city': women talk about exercise, being outdoors, and the meanings of a large urban park}.
\newblock \emph{\bibinfo{journal}{Health Place}} \textbf{\bibinfo{volume}{12}}, \bibinfo{pages}{631--643} (\bibinfo{year}{2006}).

\bibitem{Day2006-ak}
\bibinfo{author}{Day, K.}
\newblock \bibinfo{title}{Active living and social justice: Planning for physical activity in low-income, black, and {L}atino communities}.
\newblock \emph{\bibinfo{journal}{J. Am. Plann. Assoc.}} \textbf{\bibinfo{volume}{72}}, \bibinfo{pages}{88--99} (\bibinfo{year}{2006}).

\bibitem{MORALES2021}
\bibinfo{author}{ReVelle, C.}, \bibinfo{author}{Eiselt, H.} \& \bibinfo{author}{Daskin, M.}
\newblock \bibinfo{title}{Racial/ethnic disparities in household food insecurity during the covid-19 pandemic: a nationally representative study}.
\newblock \emph{\bibinfo{journal}{Racial and Ethnic Health Disparities}} \textbf{\bibinfo{volume}{8}}, \bibinfo{pages}{1300--1314} (\bibinfo{year}{2021}).
\newblock \urlprefix\url{https://link.springer.com/article/10.1007/s40615-020-00892-7}.

\bibitem{Kimani2021-xx}
\bibinfo{author}{Kimani, M.~E.}, \bibinfo{author}{Sarr, M.}, \bibinfo{author}{Cuffee, Y.}, \bibinfo{author}{Liu, C.} \& \bibinfo{author}{Webster, N.~S.}
\newblock \bibinfo{title}{Associations of {Race/Ethnicity} and food insecurity with {COVID-19} infection rates across {US} counties}.
\newblock \emph{\bibinfo{journal}{JAMA Netw Open}} \textbf{\bibinfo{volume}{4}}, \bibinfo{pages}{e2112852} (\bibinfo{year}{2021}).

\bibitem{Sharma2020-ol}
\bibinfo{author}{Sharma, S.~V.} \emph{et~al.}
\newblock \bibinfo{title}{Social determinants of {Health-Related} needs during {COVID-19} among {Low-Income} households with children}.
\newblock \emph{\bibinfo{journal}{Prev. Chronic Dis.}} \textbf{\bibinfo{volume}{17}}, \bibinfo{pages}{E119} (\bibinfo{year}{2020}).

\bibitem{pmed1}
\bibinfo{author}{Hakimi, S.}
\newblock \bibinfo{title}{Optimum locations of switching centers and the absolute centers and medians of a graph}.
\newblock \emph{\bibinfo{journal}{Operations Research}} \textbf{\bibinfo{volume}{12}}, \bibinfo{pages}{450--459} (\bibinfo{year}{1964}).

\bibitem{pmed2}
\bibinfo{author}{Hakimi, S.}
\newblock \bibinfo{title}{Optimum distribution of switching centers in a communication network and some related graph theoretic problems}.
\newblock \emph{\bibinfo{journal}{Operations Research}} \textbf{\bibinfo{volume}{13}}, \bibinfo{pages}{462--475} (\bibinfo{year}{1965}).

\bibitem{covering}
\bibinfo{author}{CHURCH, R.} \& \bibinfo{author}{VELLE, C.~R.}
\newblock \bibinfo{title}{The maximal covering location problem}.
\newblock \emph{\bibinfo{journal}{Papers of the Regional Science Association}} \textbf{\bibinfo{volume}{22}}, \bibinfo{pages}{101--118} (\bibinfo{year}{1974}).

\bibitem{Karsu2015-cb}
\bibinfo{author}{Karsu, {\"O}.} \& \bibinfo{author}{Morton, A.}
\newblock \bibinfo{title}{Inequity averse optimization in operational research}.
\newblock \emph{\bibinfo{journal}{Eur. J. Oper. Res.}} \textbf{\bibinfo{volume}{245}}, \bibinfo{pages}{343--359} (\bibinfo{year}{2015}).

\bibitem{Barbati2018-nd}
\bibinfo{author}{Barbati, M.} \& \bibinfo{author}{Bruno, G.}
\newblock \bibinfo{title}{Exploring similarities in discrete facility location models with equality measures} (\bibinfo{year}{2018}).

\bibitem{census2022}
\bibinfo{author}{{United States Census Bureau}}.
\newblock \bibinfo{title}{2020 united states census}.
\newblock \urlprefix\url{https://census.gov/programs-surveys/decennial-census/decade/2020/2020-census-main.html}.

\bibitem{Penchansky1981-qh}
\bibinfo{author}{Penchansky, R.} \& \bibinfo{author}{Thomas, J.~W.}
\newblock \bibinfo{title}{The concept of access: definition and relationship to consumer satisfaction}.
\newblock \emph{\bibinfo{journal}{Med. Care}} \textbf{\bibinfo{volume}{19}}, \bibinfo{pages}{127--140} (\bibinfo{year}{1981}).

\bibitem{Saurman2016-gj}
\bibinfo{author}{Saurman, E.}
\newblock \bibinfo{title}{Improving access: modifying penchansky and thomas's theory of access}.
\newblock \emph{\bibinfo{journal}{J. Health Serv. Res. Policy}} \textbf{\bibinfo{volume}{21}}, \bibinfo{pages}{36--39} (\bibinfo{year}{2016}).

\bibitem{Logan2019-fr}
\bibinfo{author}{Logan, T.~M.} \emph{et~al.}
\newblock \bibinfo{title}{Evaluating urban accessibility: leveraging open-source data and analytics to overcome existing limitations}.
\newblock \emph{\bibinfo{journal}{Environment and Planning B: Urban Analytics and City Science}} \textbf{\bibinfo{volume}{46}}, \bibinfo{pages}{897--913} (\bibinfo{year}{2019}).

\bibitem{osrm}
\bibinfo{author}{{Project OSRM}}.
\newblock \bibinfo{title}{{Open Source Routing Machine (OSRM)}}.
\newblock \urlprefix\url{https://project-osrm.org/}.

\bibitem{Koschinsky2016-ao}
\bibinfo{author}{Koschinsky, J.}, \bibinfo{author}{Talen, E.}, \bibinfo{author}{Alfonzo, M.} \& \bibinfo{author}{Lee, S.}
\newblock \bibinfo{title}{How walkable is walker's paradise?}
\newblock \emph{\bibinfo{journal}{Environ. Plann. B Plann. Des.}} \bibinfo{pages}{0265813515625641} (\bibinfo{year}{2016}).

\bibitem{nhgis}
\bibinfo{author}{Manson, S.} \emph{et~al.}
\newblock \bibinfo{title}{{IPUMS National Historical Geographic Information System: Version 18.0 [dataset]}}.
\newblock \bibinfo{howpublished}{\url{http://doi.org/10.18128/D050.V18.0}} (\bibinfo{year}{2023}).
\newblock \bibinfo{note}{Minneapolis, MN: IPUMS}.

\bibitem{MARSH19941}
\bibinfo{author}{Marsh, M.~T.} \& \bibinfo{author}{Schilling, D.~A.}
\newblock \bibinfo{title}{Equity measurement in facility location analysis: A review and framework}.
\newblock \emph{\bibinfo{journal}{European Journal of Operational Research}} \textbf{\bibinfo{volume}{74}}, \bibinfo{pages}{1--17} (\bibinfo{year}{1994}).
\newblock \urlprefix\url{https://www.sciencedirect.com/science/article/pii/0377221794902003}.

\bibitem{Smith2013-md}
\bibinfo{author}{Smith, H.~K.}, \bibinfo{author}{Harper, P.~R.} \& \bibinfo{author}{Potts, C.~N.}
\newblock \bibinfo{title}{Bicriteria efficiency/equity hierarchical location models for public service application} (\bibinfo{year}{2013}).

\bibitem{Batta2014-uj}
\bibinfo{author}{Batta, R.}, \bibinfo{author}{Lejeune, M.} \& \bibinfo{author}{Prasad, S.}
\newblock \bibinfo{title}{Public facility location using dispersion, population, and equity criteria}.
\newblock \emph{\bibinfo{journal}{Eur. J. Oper. Res.}} \textbf{\bibinfo{volume}{234}}, \bibinfo{pages}{819--829} (\bibinfo{year}{2014}).

\bibitem{Xu2023-iy}
\bibinfo{author}{Xu, J.}, \bibinfo{author}{Murray, A.~T.}, \bibinfo{author}{Church, R.~L.} \& \bibinfo{author}{Wei, R.}
\newblock \bibinfo{title}{Service allocation equity in location coverage analytics}.
\newblock \emph{\bibinfo{journal}{Eur. J. Oper. Res.}} \textbf{\bibinfo{volume}{305}}, \bibinfo{pages}{21--37} (\bibinfo{year}{2023}).

\bibitem{hart2011pyomo}
\bibinfo{author}{Hart, W.~E.}, \bibinfo{author}{Watson, J.-P.} \& \bibinfo{author}{Woodruff, D.~L.}
\newblock \bibinfo{title}{Pyomo: modeling and solving mathematical programs in python}.
\newblock \emph{\bibinfo{journal}{Mathematical Programming Computation}} \textbf{\bibinfo{volume}{3}}, \bibinfo{pages}{219--260} (\bibinfo{year}{2011}).

\bibitem{bynum2021pyomo}
\bibinfo{author}{Bynum, M.~L.} \emph{et~al.}
\newblock \emph{\bibinfo{title}{Pyomo--optimization modeling in python}}, vol.~\bibinfo{volume}{67} (\bibinfo{publisher}{{Springer Science \& Business Media}}, \bibinfo{year}{2021}), \bibinfo{edition}{{Third}} edn.

\bibitem{gurobi}
\bibinfo{author}{{Gurobi Optimization, LLC}}.
\newblock \bibinfo{title}{{Gurobi Optimizer Reference Manual}} (\bibinfo{year}{2023}).

\bibitem{Moreno2021-vs}
\bibinfo{author}{Moreno, C.}, \bibinfo{author}{Allam, Z.}, \bibinfo{author}{Chabaud, D.}, \bibinfo{author}{Gall, C.} \& \bibinfo{author}{Pratlong, F.}
\newblock \bibinfo{title}{Introducing the ``15-minute city'': Sustainability, resilience and place identity in future {Post-Pandemic} cities}.
\newblock \emph{\bibinfo{journal}{Smart Cities}} \textbf{\bibinfo{volume}{4}}, \bibinfo{pages}{93--111} (\bibinfo{year}{2021}).

\bibitem{Capasso_Da_Silva2019-ge}
\bibinfo{author}{Capasso Da~Silva, D.}, \bibinfo{author}{King, D.~A.} \& \bibinfo{author}{Lemar, S.}
\newblock \bibinfo{title}{Accessibility in practice: 20-minute city as a sustainability planning goal}.
\newblock \emph{\bibinfo{journal}{Sustain. Sci. Pract. Policy}} \textbf{\bibinfo{volume}{12}}, \bibinfo{pages}{129} (\bibinfo{year}{2019}).

\bibitem{McNeil2011-zv}
\bibinfo{author}{McNeil, N.}
\newblock \bibinfo{title}{Bikeability and the 20-min neighborhood: How infrastructure and destinations influence bicycle accessibility}.
\newblock \emph{\bibinfo{journal}{Transp. Res. Rec.}} \textbf{\bibinfo{volume}{2247}}, \bibinfo{pages}{53--63} (\bibinfo{year}{2011}).

\bibitem{Schaap2016-ov}
\bibinfo{author}{Schaap, N.}, \bibinfo{author}{Harms, L.}, \bibinfo{author}{Kansen, M.} \& \bibinfo{author}{W{\"u}st, H.}
\newblock \bibinfo{title}{Cycling and walking: the grease in our mobility chain}.
\newblock \bibinfo{type}{Tech. Rep.} \bibinfo{number}{KiM-16-A03}, \bibinfo{institution}{Dutch Ministry of Infrastructure and the Environment} (\bibinfo{year}{2016}).

\bibitem{Frumkin2004-yi}
\bibinfo{author}{Frumkin, H.}, \bibinfo{author}{Frank, L.} \& \bibinfo{author}{Jackson, R.~J.}
\newblock \emph{\bibinfo{title}{Urban Sprawl and Public Health: Designing, Planning, and Building for Healthy Communities}} (\bibinfo{publisher}{Island Press}, \bibinfo{year}{2004}).

\bibitem{Gusdorf2007-ms}
\bibinfo{author}{Gusdorf, F.} \& \bibinfo{author}{Hallegatte, S.}
\newblock \bibinfo{title}{Compact or spread-out cities: Urban planning, taxation, and the vulnerability to transportation shocks}.
\newblock \emph{\bibinfo{journal}{Energy Policy}} \textbf{\bibinfo{volume}{35}}, \bibinfo{pages}{4826--4838} (\bibinfo{year}{2007}).

\bibitem{Patel2020-yq}
\bibinfo{author}{Patel, J.~A.} \emph{et~al.}
\newblock \bibinfo{title}{Poverty, inequality and {COVID-19}: the forgotten vulnerable}.
\newblock \emph{\bibinfo{journal}{Public Health}} \textbf{\bibinfo{volume}{183}}, \bibinfo{pages}{110--111} (\bibinfo{year}{2020}).

\bibitem{Power2020-yb}
\bibinfo{author}{Power, M.}, \bibinfo{author}{Doherty, B.}, \bibinfo{author}{Pybus, K.~J.} \& \bibinfo{author}{Pickett, K.~E.}
\newblock \bibinfo{title}{How {COVID-19} has exposed inequalities in the {UK} food system: the case of {UK} food and poverty}.
\newblock \emph{\bibinfo{journal}{Emerald Open Research}} \textbf{\bibinfo{volume}{1}} (\bibinfo{year}{2020}).

\bibitem{Rose2011-gy}
\bibinfo{author}{Rose, D.}, \bibinfo{author}{Bodor, J.~N.}, \bibinfo{author}{Rice, J.~C.}, \bibinfo{author}{Swalm, C.~M.} \& \bibinfo{author}{Hutchinson, P.~L.}
\newblock \bibinfo{title}{The effects of hurricane katrina on food access disparities in new orleans}.
\newblock \emph{\bibinfo{journal}{Am. J. Public Health}} \textbf{\bibinfo{volume}{101}}, \bibinfo{pages}{482--484} (\bibinfo{year}{2011}).

\end{thebibliography}
